# Central limit theorems for eigenvalues in a spiked population model

## Zhidong Bai[1,a] and Jian-feng Yao[2,b]

[a] KLASMOE and School of Mathematics and Statistics, Northeast Normal University, 5268 People's Road, 130024
Changchun, China and Department of Statistics and Applied Probability, National University of Singapore, 10, Kent Ridge
Crescent, Singapore 119260. E-mail: stabaizd@nus.edu.sg
[b] IRMAR and Université de Rennes 1, Campus de Beaulieu, 35042 Rennes Cedex, France.
E-mail: jian-feng.yao@univ-rennes1.fr

**Abstract.** In a spiked population model, the population covariance matrix has all its eigenvalues equal to units except for a few fixed eigenvalues (spikes). This model is proposed by Johnstone to cope with empirical findings on various data sets. The question is to quantify the effect of the perturbation caused by the spike eigenvalues. A recent work by Baik and Silverstein establishes the almost sure limits of the extreme sample eigenvalues associated to the spike eigenvalues when the population and the sample sizes become large. This paper establishes the limiting distributions of these extreme sample eigenvalues. As another important result of the paper, we provide a central limit theorem on random sesquilinear forms.

**Résumé.** Dans un modèle de variances hétérogènes, les valeurs propres de la matrice de covariance des variables sont toutes égales à l'unité sauf un faible nombre d'entre elles. Ce modèle a été introduit par Johnstone comme une explication possible de la structure des valeurs propres de la matrice de covariance empirique constatée sur plusieurs ensembles de données réelles. Une question importante est de quantifier la perturbation causée par ces valeurs propres différentes de l'unité. Un travail récent de Baik et Silverstein établit la limite presque sûre des valeurs propres empiriques extrêmes lorsque le nombre de variables tend vers l'infini proportionnellement à la taille de l'échantillon. Ce travail établit un théorème limite central pour ces valeurs propres empiriques extrêmes. Il est basé sur un nouveau théorème limite central pour les formes sesquilinéaires aléatoires.

## 1. Introduction

It is well known that the empirical spectral distribution (ESD) of a large sample covariance matrix converges to the family of Marčenko–Pastur laws under fairly general condition on the sample variables [3, 10]. On the other hand, the study of the largest or smallest eigenvalues is more complex. In a variety of situations, the almost sure limits of these extreme eigenvalues are proved to coincide with the boundaries of the support of the limiting distribution. As an example, when the sample vectors have independent coordinates and unit variances and assuming that the ratio $p/n$ of the population size $p$ over the sample size $n$ tends to a positive limit $y \in (0, 1)$, then the limiting distribution is the classical Marčenko–Pastur law $F_y(\mathrm{d}x)$

$$F_y(\mathrm{d}x) = \frac{1}{2\pi xy}\sqrt{(x - a_y)(b_y - x)}\,\mathrm{d}x, \quad a_y \le x \le b_y, \tag{1.1}$$

where $a_y = (1 - \sqrt{y})^2$, and $b_y = (1 + \sqrt{y})^2$. Moreover, the smallest and the largest eigenvalue converge almost surely to the boundary $a_y$ and $b_y$, respectively.

---

[1]The research of this author was supported by CNSF grant 10571020 and NUS grant R-155-000-061-112.

[2]Research was (partially) completed while J.-F. Yao was visiting the Institute for Mathematical Sciences, National University of Singapore in 2006.

Recent empirical data analysis from fields like wireless communication engineering, speech recognition or gene expression experiments suggest that frequently, some extreme eigenvalues of sample covariance matrices are well-separated from the rest. For instance, see Figs 1 and 2 in [9] which display the sample eigenvalues of the functional data consisting of a speech data set of 162 instances of a phoneme "dcl" spoken by males calculated at 256 points. As a way for possible explanation of this phenomenon, this author proposes a *spiked population model* where all eigenvalues of the population covariance matrix are equal to one except a fixed and relatively small number among them (*spikes*). Clearly, a spiked population model can be considered as a small perturbation of the so-called *null case* where all the eigenvalues of the population covariance matrix are unit. It then raises the question how such a small perturbation affects the limits of the extreme eigenvalues of the sample covariance matrix as compared to the *null case*.

The behavior of the largest eigenvalue in the case of complex Gaussian variables has been recently studied in [7]. These authors prove a transition phenomenon: the weak limit as well as the scaling of the largest eigenvalue is different according to whether the largest spike eigenvalue is larger, equal or less than the critical value $1 + \sqrt{y}$. In [6], the authors consider the spiked population model with general random variables: complex or real and not necessarily Gaussian. For the almost sure limits of the extreme sample eigenvalues, they also find that these limits depend on the critical values $1 + \sqrt{y}$ and $1 - \sqrt{y}$ from above and below, respectively. For example, if there are $M$ eigenvalues in the population covariance matrix larger than $1 + \sqrt{y}$, then the $M$ largest eigenvalues from the sample covariance matrix will (almost surely) have their limits above the right edge $b_y$ of the limiting Marčenko–Pastur law. Analogous results are also proposed for the case $y > 1$ and $y = 1$.

An important question here is to find the limiting distributions of these extreme eigenvalues. As mentioned above, the results are proposed in [7] for the largest eigenvalue and the Gaussian complex case. In this perspective, assuming that the population vector is real Gaussian with a diagonal covariance matrix and that the $M$ spike eigenvalues are all simple, [12] found that each of the $M$ largest sample eigenvalues has a Gaussian limiting distribution.

In this paper, we follow the general set-up of [6]. Assuming $y \in (0, 1)$ and general population variables, we will establish central limit theorems for the largest as well as for the smallest sample eigenvalues associated to spike eigenvalues outside the interval $[1 - \sqrt{y}, 1 + \sqrt{y}]$. Furthermore, we prove that the limiting distribution of such sample extreme eigenvalues is Gaussian only if the corresponding spike population eigenvalue is simple. Otherwise, if a spiked eigenvalue is multiple, say of index $k$, then there will be $k$ packed-consecutive sample eigenvalues $\lambda_{n,1}, \ldots, \lambda_{n,k}$ which converge jointly to the distribution of a $k \times k$ symmetric (or Hermitian) Gaussian random matrix. Consequently in this case, the limiting distribution of a single $\lambda_{n,j}$ is generally non Gaussian.

The main tools of our analysis are borrowed from the random matrix theory on one hand. For general background of this theory, we refer to the book [11] and a modern review by Bai [3]. On the other hand, we introduce in this paper another important tool, namely a CLT for random sesquilinear forms which should have its own interests. This CLT, independent from the rest of the paper, is presented in the last section (Section 7).

The remaining sections of the paper are organized as follows. First in Section 2, we introduce the spiked population model and recall known results on the almost sure limits of extreme sample eigenvalues. The main result of the paper, namely a general CLT for extreme sample eigenvalues, Theorem 3.1, is then introduced in Section 3. To provide a better account of this CLT, Section 4 develops in details several meaningful examples. Several sets of numerical computations are also conducted to give concrete illustration of the main result. In particular, we recover a CLT given in [12] as a special instance. In Section 5, we discuss some extensions of these results to the case where spiked eigenvalues are inside the gaps located in the center of the spectrum of the population covariance matrix. Finally, Section 6 collects the proofs of the presented results based on a CLT for random sesquilinear forms which is itself introduced and proved in Section 7.

## 2. Spiked population model and convergence of extreme eigenvalues

We consider a zero-mean, complex-valued random vector $x = (\xi^{\mathrm{T}}, \eta^{\mathrm{T}})^{\mathrm{T}}$ where $\xi = (\xi(1), \ldots, \xi(M))^{\mathrm{T}}$, $\eta = (\eta(1), \ldots, \eta(p))^{\mathrm{T}}$ are independent, of dimension $M$ and $p$, respectively. Moreover, we assume that $\mathbb{E}[\|x\|^4] < \infty$ and the coordinates of $\eta$ are independent and identically distributed with unit variance. The *population covariance matrix* of the vector $x$ is therefore

$$V = \mathrm{cov}(x) = \begin{pmatrix} \varSigma & 0 \\ 0 & I_p \end{pmatrix}.$$

We consider the following spiked population model by assuming that $\varSigma$ has $K$ non null and non unit eigenvalues $\alpha_1, \ldots, \alpha_K$ with respective multiplicity $n_1, \ldots, n_K$ $(n_1 + \cdots + n_K = M)$. Therefore, the eigenvalues of the population covariance matrix $V$ are units except the $(\alpha_j)$, called *spike eigenvalues*.

Let $x_i = (\xi_i^{\mathrm{T}}, \eta_i^{\mathrm{T}})^{\mathrm{T}}$ be $n$ copies i.i.d. of $x$. The *sample covariance matrix* is

$$S_n = \frac{1}{n} \sum_{i=1}^n x_i x_i^*,$$

which can be rewritten as

$$S_n = \begin{pmatrix} S_{11} & S_{12} \\ S_{21} & S_{22} \end{pmatrix} = \begin{pmatrix} X_1 X_1^* & X_1 X_2^* \\ X_2 X_1^* & X_2 X_2^* \end{pmatrix} = \frac{1}{n} \begin{pmatrix} \sum \xi_i \xi_i^* & \sum \xi_i \eta_i^* \\ \sum \eta_i \xi_i^* & \sum \eta_i \eta_i^* \end{pmatrix}, \tag{2.1}$$

with

$$X_1 = \frac{1}{\sqrt{n}} (\xi_1, \ldots, \xi_n)_{M \times n} = \frac{1}{\sqrt{n}} \xi_{1:n}, \qquad X_2 = \frac{1}{\sqrt{n}} (\eta_1, \ldots, \eta_n)_{p \times n} = \frac{1}{\sqrt{n}} \eta_{1:n}.$$

It is assumed in the sequel that $M$ is fixed, and $p$ and $n$ are related so that when $n \to \infty$, $p/n \to y \in (0, 1)$. The ESD of $S_n$, as well as the one of $S_{22}$, converges to the Marčenko–Pastur distribution $F_y(\mathrm{d}x)$ given in (1.1). As explained in the Introduction, a central question is to quantify the effect caused by the small number of spiked eigenvalues on the asymptotic of the extreme sample eigenvalues.

As a first general answer to this question, Baik and Silverstein [6] completely determines the almost sure limits of largest and smallest sample eigenvalues. More precisely, assume that among the $M$ eigenvalues of $\varSigma$, there are exactly $M_b$ greater than $1 + \sqrt{y}$ and $M_a$ smaller than $1 - \sqrt{y}$:

$$\alpha_1 > \cdots > \alpha_{M_b} > 1 + \sqrt{y}, \qquad \alpha_M < \cdots < \alpha_{M-M_b+1} < 1 - \sqrt{y}, \tag{2.2}$$

and $1 - \sqrt{y} \le \alpha_k \le 1 + \sqrt{y}$ for the other $\alpha_k$'s. Moreover, for $\alpha \ne 1$, we define the function

$$\lambda = \phi(\alpha) = \alpha + \frac{y\alpha}{\alpha - 1}. \tag{2.3}$$

As $y < 1$, we have $p \le n$ for large $n$. Let

$$\lambda_{n,1} \ge \lambda_{n,2} \ge \cdots \ge \lambda_{n,p}$$

be the eigenvalues of the sample covariance matrix $S_n$. Let $s_i = n_1 + \cdots + n_i$ for $1 \le i \le M_b$ and $t_j = n_M + \cdots + n_j$ for $1 \le j \le M_a$ (by convention $s_0 = t_0 = 0$).

Therefore, Baik and Silverstein [6] proves that for each $k \in \{1, \ldots, M_b\}$ and $s_{k-1} < j \le s_k$ (largest eigenvalues) or $k \in \{1, \ldots, M_a\}$ and $p - t_k < j \le p - t_{k-1}$ (smallest eigenvalues),

$$\lambda_{n,j} \to \phi(\alpha_k) = \alpha_k + \frac{y\alpha_k}{\alpha_k - 1}, \quad \text{almost surely.} \tag{2.4}$$

In other words, if a spike eigenvalue $\alpha_k$ lies outside the interval $[1 - \sqrt{y}, 1 + \sqrt{y}]$ and has multiplicity $n_k$, then $\phi(\alpha_k)$ is the limit of $n_k$ packed sample eigenvalue $\{\lambda_{n,j}, j \in J_k\}$. Here we have denoted by $J_k$ the corresponding set of indexes: $J_k = \{s_{k-1} + 1, \ldots, s_k\}$ for $\alpha_k > 1 + \sqrt{y}$ and $J_k = \{p - t_k + 1, \ldots, p - t_{k-1}\}$ for $\alpha_k < 1 - \sqrt{y}$.

## 3. Main results

The aim of this paper is to derive a CLT for the $n_k$-packed sample eigenvalues

$$\sqrt{n}[\lambda_{n,j} - \phi(\alpha_k)], \quad j \in J_k,$$

where $\alpha_k \notin [1 - \sqrt{y}, 1 + \sqrt{y}]$ is some fixed spike eigenvalue of multiplicity $n_k$. The statement of the main result of the paper, Theorem 3.1, needs several intermediate notations and results.

### 3.1. Determinant equation and a random sesquilinear form

By definition, each $\lambda_{n,j}$ solves the equation

$$0 = |\lambda I - S_n| = |\lambda I - S_{22}||\lambda I - K_n(\lambda)|, \tag{3.1}$$

where

$$K_n(\lambda) = S_{11} + S_{12}(\lambda I - S_{22})^{-1}S_{21}. \tag{3.2}$$

As when $n \to \infty$, with probability 1, the limit $\lambda_{n,j} \to \phi(\alpha_k) \notin [a_y, b_y]$ and the eigenvalues of $S_{22}$ go inside the interval $[a_y, b_y]$, the probability of the event $Q_n$

$$Q_n = \{\lambda_{n,j} \notin [a_y, b_y]\} \cap \{\text{spectrum of } S_{22} \subset [a_y, b_y]\}$$

tends to 1. Conditional on this event, the $(\lambda_{n,j})$'s then solve the *determinant equation*

$$|\lambda I - K_n(\lambda)| = 0. \tag{3.3}$$

Therefore without loss of generality, we can assume that $\lambda_{n,j} \notin [a_y, b_y]$ and they are solutions of this equation.

Furthermore, let

$$A_n = (a_{ij}) = A_n(\lambda) = X_2^*(\lambda I - X_2 X_2^*)^{-1}X_2, \quad \lambda \notin [a_y, b_y]. \tag{3.4}$$

Lemma 6.1 detailed in Section 6.1 establishes the convergence of several statistics of the matrix $A_n$. In particular, $n^{-1}\operatorname{tr} A_n$, $n^{-1}\operatorname{tr} A_n A_n^*$ and $n^{-1}\sum_{i=1}^{n} a_{ii}^2$ converges in probability to $ym_1(\lambda)$, $ym_2(\lambda)$ and $(y[1 + m_1(\lambda)]/\{\lambda - y[1 + m_1(\lambda)]\})^2$, respectively. Here, the $m_j(\lambda)$ are some specific transforms of the Marčenko–Pastur law $F_y$ (see Section 6.1 for more details).

Therefore, the random form $K_n$ in (3.2) can be decomposed as follows

$$\begin{aligned} K_n(\lambda) &= S_{11} + X_1 A_n X_1^* = \frac{1}{n}\xi_{1:n}(I + A_n)\xi_{1:n}^* \\ &= \frac{1}{n}\{\xi_{1:n}(I + A_n)\xi_{1:n}^* - \Sigma\operatorname{tr}(I + A_n)\} + \frac{1}{n}\Sigma\operatorname{tr}(I + A_n) \\ &= \frac{1}{\sqrt{n}}R_n + [1 + ym_1(\lambda)]\Sigma + o_P\left(\frac{1}{\sqrt{n}}\right), \end{aligned} \tag{3.5}$$

with

$$R_n = R_n(\lambda) = \frac{1}{\sqrt{n}} \{ \xi_{1:n}(I + A_n)\xi_{1:n}^* - \Sigma \operatorname{tr}(I + A_n) \}. \tag{3.6}$$

In the last derivation, we have used the fact

$$\frac{1}{n} \operatorname{tr}(I + A_n) = 1 + y m_1(\lambda) + \mathrm{o}_P\left(\frac{1}{\sqrt{n}}\right),$$

which follows from a CLT for $\operatorname{tr}(A_n)$ (see [4]).

### 3.2. Limit distribution of the random matrices $\{R_n(\lambda)\}$

The next step is to find the limit distribution of the sequence of random matrices $\{R_n(\lambda)\}$. The situation is different for the real and complex cases. Define the constants

$$\theta = 1 + 2y m_1(\lambda) + y m_2(\lambda), \tag{3.7}$$

$$\omega = 1 + 2y m_1(\lambda) + \left( \frac{y[1 + m_1(\lambda)]}{\lambda - y[1 + m_1(\lambda)]} \right)^2. \tag{3.8}$$

**Proposition 3.1 (Limiting distribution of $R_n(\lambda)$: real variables case).** *Assume that the variables $\xi$ and $\eta$ are real-valued. Then, the random matrix $R_n$ converges weakly to a symmetric random matrix $R = (R_{ij})$ with zero-mean Gaussian entries having the following covariance function: for $1 \le i \le j \le M$, $1 \le i' \le j' \le M$*

$$\operatorname{cov}(R_{ij}, R_{i'j'}) = \omega \{ \mathbb{E}[\xi(i)\xi(j)\xi(i')\xi(j')] - \Sigma_{ij}\Sigma_{i'j'} \} + (\theta - \omega)\{ \mathbb{E}[\xi(i)\xi(j')]\mathbb{E}[\xi(i')\xi(j)] \}$$
$$+ (\theta - \omega)\{ \mathbb{E}[\xi(i)\xi(i')]\mathbb{E}[\xi(j)\xi(j')] \}. \tag{3.9}$$

Note that in particular, the following formula holds for the variances

$$\operatorname{var}(R_{ij}) = \theta(\Sigma_{ii}\Sigma_{jj} + \Sigma_{ij}^2) + \omega\{ \mathbb{E}[\xi^2(i)\xi^2(j)] - 2\Sigma_{ij}^2 - \Sigma_{ii}\Sigma_{jj} \}. \tag{3.10}$$

In case of a diagonal element $R_{ii}$, this expression simplifies to

$$\operatorname{var}(R_{ii}) = [2\theta + \beta_i \omega]\Sigma_{ii}^2, \quad \text{with } \beta_i = \frac{\mathbb{E}[\xi(i)^4]}{\Sigma_{ii}^2} - 3. \tag{3.11}$$

If moreover, $\xi(i)$ is Gaussian, $\beta_i = 0$.

**Remark 1.** *If the coordinates $\{\xi(i)\}$ of $\xi$ are independent, then the limiting covariance matrix in (3.9) is diagonal: the limiting Gaussian matrix is made with independent entries. Their variances simplify to (3.11) and*

$$\operatorname{var}(R_{ij}) = \theta \Sigma_{ii}\Sigma_{jj}, \quad i < j. \tag{3.12}$$

**Proposition 3.2 (Limiting distribution of $R_n(\lambda)$: complex variables case).** *Assume the general case with complex-valued variables $\xi$ and $\eta$ and that the following limit exists*

$$m_4(\lambda) = \lim_n \frac{1}{n} \operatorname{tr} A_n A_n^{\mathrm{T}}, \quad \lambda \notin [a_y, b_y]. \tag{3.13}$$

*Then, the random matrix $R_n$ converges weakly to a zero-mean Hermitian random matrix $R = (R_{ij})$. Moreover, the joint distribution of the real and imaginary parts of the upper-triangular bloc $\{R_{ij}, 1 \le i \le j \le M\}$ is a $2K$-dimensional Gaussian vector with covariance matrix*

$$\Gamma = \begin{pmatrix} \Gamma_{11} & \Gamma_{12} \\ \Gamma_{21} & \Gamma_{22} \end{pmatrix}, \tag{3.14}$$

*where*

$$\Gamma_{11} = \frac{1}{4} \sum_{j=1}^{3} \{2\Re(B_j) + B_{ja} + B_{jb}\},$$

$$\Gamma_{22} = \frac{1}{4} \sum_{j=1}^{3} \{-2\Re(B_j) + B_{ja} + B_{jb}\},$$

$$\Gamma_{12} = \frac{1}{2} \sum_{j=1}^{3} \Im(B_j),$$

*and for $1 \le i \le j \le M$ and $1 \le i' \le j' \le M$,*

$$B_1(ij, i'j') = \omega(\mathbb{E}[\xi(i)\bar{\xi}(j)\xi(i')\bar{\xi}(j')] - \Sigma_{ij}\Sigma_{i'j'}),$$

$$B_2(ij, i'j') = (\theta - \omega)\Sigma_{ij'}\Sigma_{i'j},$$

$$B_3(ij, i'j') = (\tau - \omega)(\mathbb{E}[\xi(i)\xi(i')]\mathbb{E}[\bar{\xi}(j)\bar{\xi}(j')]),$$

$$B_{1a}(ij, i'j') = \omega(\mathbb{E}[|\xi(i)\xi(i')|^2] - \Sigma_{ii}\Sigma_{i'i'}),$$

$$B_{1b}(ij, i'j') = \omega(\mathbb{E}[|\xi(j)\xi(j')|^2] - \Sigma_{jj}\Sigma_{j'j'}),$$

$$B_{2a}(ij, i'j') = (\theta - \omega)|\Sigma_{ii'}|^2,$$

$$B_{2b}(ij, i'j') = (\theta - \omega)|\Sigma_{jj'}|^2,$$

$$B_{3a}(ij, i'j') = (\tau - \omega)|\mathbb{E}[\xi(i)\xi(i')]|^2,$$

$$B_{3b}(ij, i'j') = (\tau - \omega)|\mathbb{E}[\xi(j)\xi(j')]|^2.$$

*Here, the constant $\tau$ equals*

$$\tau = \lim_n \frac{1}{n} \operatorname{tr}(I + A_n)(I + A_n)^{\mathrm{T}} = 1 + 2ym_1(\lambda) + m_4(\lambda). \tag{3.15}$$

The limiting covariance matrix $\Gamma$ has a complicated expression. However, the variance of a diagonal element $R_{ii}$ has a much simpler expression if moreover, $\mathbb{E}[\xi^2(i)] = 0$ for all $1 \le i \le M$,

$$\operatorname{var}(R_{ii}) = [\theta + \beta_i'\omega]\Sigma_{ii}^2, \quad \text{with } \beta_i' = \frac{\mathbb{E}[\xi(i)^4]}{\Sigma_{ii}^2} - 2. \tag{3.16}$$

In particular, if $\xi(i)$ is Gaussian, $\beta_i' = 0$.

### 3.3. CLT for extreme eigenvalues

In order to introduce the main result of the paper, let the spectral decomposition of $\Sigma$,

$$\Sigma = U \begin{pmatrix} \alpha_1 I_{n_1} & \cdots & 0 \\ 0 & \ddots & 0 \\ \cdots & 0 & \alpha_K I_{n_K} \end{pmatrix} U^*, \tag{3.17}$$

where $U$ is an unitary matrix. Following Section 2, for each spiked eigenvalue $\alpha_k \notin [1 - \sqrt{y}, 1 + \sqrt{y}]$, let $\{\lambda_{n,j}, j \in J_k\}$ be the $n_k$ packed eigenvalues of the sample covariance matrix which all tend almost surely to $\lambda_k = \phi(\alpha_k)$. Let $R(\lambda_k)$ be the Gaussian matrix limit of the sequence of matrices of random forms $[R_n(\lambda_k)]_n$ given in Proposition 3.1 (real variables case) and Proposition 3.2 (complex variables case), respectively. Let

$$\widetilde{R}(\lambda_k) = U^* R(\lambda_k) U. \tag{3.18}$$

**Theorem 3.1.** *For each spike eigenvalue $\alpha_k \notin [1 - \sqrt{y}, 1 + \sqrt{y}]$, the $n_k$-dimensional real vector*

$$\sqrt{n}\{\lambda_{n,j} - \lambda_k, j \in J_k\},$$

*converges weakly to the distribution of the $n_k$ eigenvalues of the Gaussian random matrix*

$$\frac{1}{1 + y m_3(\lambda_k)\alpha_k} \widetilde{R}_{kk}(\lambda_k),$$

*where $\widetilde{R}_{kk}(\lambda_k)$ is the $k$th diagonal bloc of $\widetilde{R}(\lambda_k)$ corresponding to the indexes $\{u, v \in J_k\}$.*

One striking fact from this theorem is that the limiting distribution of such $n_k$ packed sample extreme eigenvalues are generally *non-Gaussian* and asymptotically dependent. Indeed, the limiting distribution of a single sample extreme eigenvalue $\lambda_{n,j}$ is Gaussian if and only if the corresponding population spike eigenvalue is simple.

## 4. Examples and numerical results

This section is devoted to describe in more details the content of Theorem 3.1 with several meaningful examples together with extended numerical computations.

### 4.1. A special Gaussian case from Paul [12]

We consider a particular situation examined in [12]. Assume that the variables are real Gaussian, $\Sigma$ diagonal whose eigenvalues are all simple. In other words, $K = M$ and $n_k = 1$ for all $1 \leq k \leq M$. Hence, $U = I_M$. Following Theorem 3.1, for any $\lambda_k = \phi(\alpha_k)$, with $\alpha_k \notin [1 \mp \sqrt{y}]$ $\sqrt{n}(\lambda_{n,k} - \lambda_k)$ converges weakly to the Gaussian variable $(1 + y m_3(\lambda_k)\alpha_k)^{-1} R(\lambda_k)_{kk}$. This variable is zero-mean. For the computation of its variance, we remark that by Eq. (3.11)

$$\operatorname{var} R(\lambda_k)_{kk} = 2\theta \alpha_k^2,$$

where

$$\theta = 1 + 2y m_1(\lambda_k) + y m_2(\lambda_k) = \frac{(\alpha_k - 1 + y)^2}{(\alpha_k - 1)^2 - y}.$$

Taking into account (6.6), we get finally, for $1 \leq k \leq M$

$$\sqrt{n}(\lambda_{n,k} - \lambda_k) \xrightarrow{\mathcal{D}} \mathcal{N}(0, \sigma_{\alpha_k}^2), \quad \sigma_{\alpha_k}^2 = \frac{2\alpha_k^2[(\alpha_k - 1)^2 - y]}{(\alpha_k - 1)^2}.$$

This coincides with Theorem 3 of [12].

### 4.2. More general Gaussian variables case

In this example, we assume that all variables are real Gaussian, and the coordinates of $\xi$ are independent. As in [6], we fix $y = 0.5$. The critical interval is then $[1 - \sqrt{y}, 1 + \sqrt{y}] = [0.293, 1.707]$ and the limiting support $[a_y, b_y] = [0.086, 2.914]$.

Consider $K = 4$ spike eigenvalues $(\alpha_1, \alpha_2, \alpha_3, \alpha_4) = (4, 3, 0.2, 0.1)$ with respective multiplicity $(n_1, n_2, n_3, n_4) = (1, 2, 2, 1)$. Let

$$\lambda_{n,1} \geq \lambda_{n,2} \geq \lambda_{n,3} \quad \text{and} \quad \lambda_{n,4} \geq \lambda_{n,5} \geq \lambda_{n,6}$$

be, respectively, the three largest and the three smallest eigenvalues of the sample covariance matrix. Let, as in Section 4.1,

$$\sigma_{\alpha_k}^2 = \frac{2\alpha_k^2[(\alpha_k - 1)^2 - y]}{(\alpha_k - 1)^2}. \tag{4.1}$$

We have $(\sigma_{\alpha_k}^2, k = 1, \dots, 4) = (30.222, 15.75, 0.0175, 0.00765)$.

Following Theorem 3.1, taking into account Section 4.1 and Proposition 3.1, we have

- for $j = 1$ and 6,

$$\delta_{n,j} = \sqrt{n}[\lambda_{n,j} - \phi(\alpha_k)] \overset{\mathcal{D}}{\Longrightarrow} \mathcal{N}(0, \sigma_{\alpha_k}^2). \tag{4.2}$$

  Here, for $j = 1$, $k = 1$, $\phi(\alpha_1) = 4.667$ and $\sigma_{\alpha_1}^2 = 30.222$; and for $j = 6$, $k = 4$, $\phi(\alpha_4) = 0.044$ and $\sigma_{\alpha_4}^2 = 0.00765$;

- for $j = (2,3)$ or $j = (4,5)$, the two-dimensional vector $\delta_{n,j} = \sqrt{n}[\lambda_{n,j} - \phi(\alpha_k)]$ converges weakly to the distribution of (ordered) eigenvalues of the random matrix

$$G = \sigma_{\alpha_k} \begin{pmatrix} W_{11} & W_{12} \\ W_{12} & W_{22} \end{pmatrix}.$$

  Here, because the initial variables $(\xi(i))$'s are Gaussian, by Eqs (3.11) and (3.12), we have $\mathrm{var}(W_{11}) = \mathrm{var}(W_{22}) = 1$, $\mathrm{var}(W_{12}) = \frac{1}{2}$, so that $(W_{ij})$ is a real Gaussian–Wigner matrix (with independent entries). Again, the variance parameter $\sigma_{\alpha_k}^2$ is defined as previously but with $k = 2$ for $j = (2,3)$ and $k = 3$ for $j = (4,5)$, respectively. Since the joint distribution of eigenvalues of a Gaussian–Wigner matrix is known (see [11]), we get the following (unordered) density for the limiting distribution of $\delta_{n,j}$:

$$g(\delta, \gamma) = \frac{1}{4\sigma_{\alpha_k}^3 \sqrt{\pi}} |\delta - \gamma| \exp\left[ -\frac{1}{2\sigma_{\alpha_k}^2} (\delta^2 + \gamma^2) \right]. \tag{4.3}$$

Experiments are conducted to compare numerically the empirical distribution of the $\delta_{n,j}$'s to their limiting value. To this end, we fix $p = 500$ and $n = 1000$. We repeat 1000 independent simulations to get 1000 replications of the six random variates $\{\delta_{n,j}, j = 1, \dots, 6\}$. Based on these replications, we compute

- a kernel density estimate for two univariate variables $\delta_{n,1}$ and $\delta_{n,6}$, denoted by $\widehat{f}_{n,1}$ $\widehat{f}_{n,6}$ respectively;
- a kernel density estimate for two bivariate variables $(\delta_{n,2}, \delta_{n,3})$ and $(\delta_{n,4}, \delta_{n,5})$, denoted by $\widehat{f}_{n,23}$ $\widehat{f}_{n,45}$ respectively.

The kernel density estimates are computed using the R software implementing an automatic bandwidth selection method from [13].

Figure 1 compare the two univariate density estimates $\widehat{f}_{n,1}$ and $\widehat{f}_{n,6}$ to their Gaussian limits (4.2). As we can see, the simulations confirm well the found formula.

To compare the bivariate density estimates $\widehat{f}_{n,23}$ and $\widehat{f}_{n,45}$ to their limiting densities given in (4.3), we choose to display their contour lines. This is done in Fig. 2 for $\widehat{f}_{n,23}$ and Fig. 3 for $\widehat{f}_{n,45}$. Again we see that the theoretical result is well confirmed.

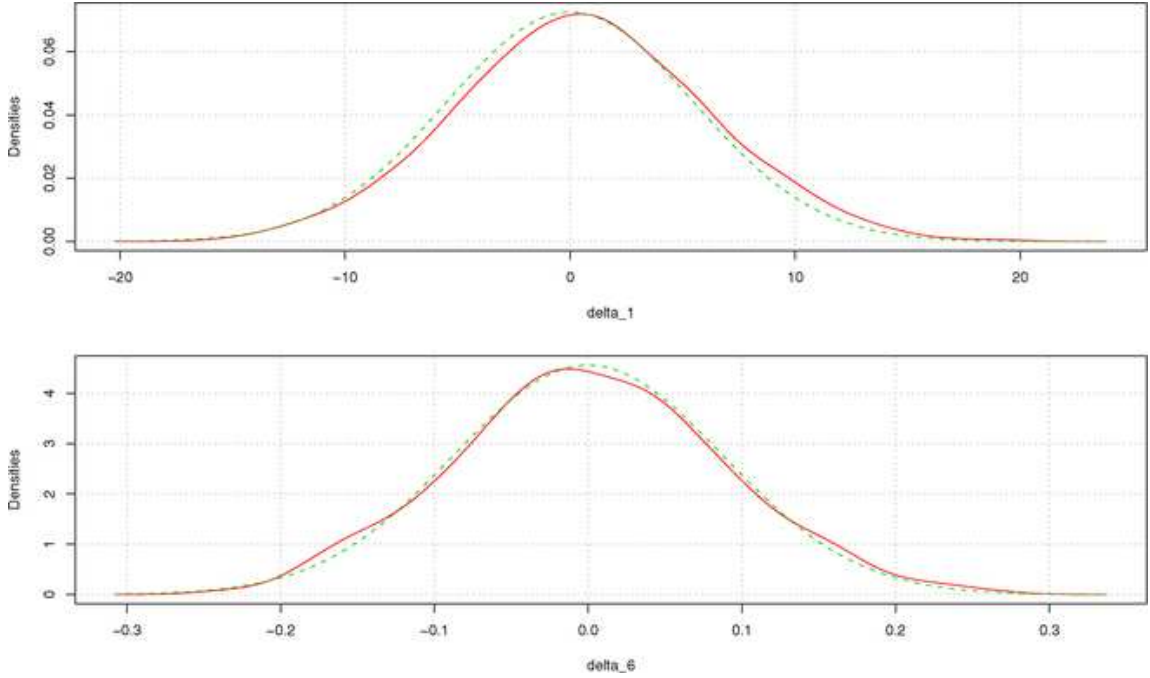

Fig. 1. Empirical density estimates (in solid lines) from the largest (top: $\widehat{f}_{n,1}$) and the smallest (bottom: $\widehat{f}_{n,6}$) sample eigenvalue from 1000 independent replications, compared to their Gaussian limits (dashed lines). Gaussian entries with $p = 500$ and $n = 1000$.

### 4.3. A binary variables case

As in the previous example, we fix $y = 0.5$ and adopt the same spike eigenvalues $(\alpha_1, \alpha_2, \alpha_3, \alpha_4) = (4, 3, 0.2, 0.1)$ with multiplicities $(n_1, n_2, n_3, n_4) = (1, 2, 2, 1)$. Let the $\sigma_{\alpha_k}^2$'s be as defined in (4.1). Again we assume that all the coordinates are independent but this time we consider binary entries. To cope with the eigenvalues, we set

$$\xi(i) = \sqrt{\alpha_k} \varepsilon_i, \qquad \eta(j) = \varepsilon_j',$$

where $(\varepsilon_i)$ and $(\varepsilon_j')$ are two independent sequences of i.i.d. binary variables taking values $\{+1, -1\}$ with equiprobability. We remark that $\mathbb{E}\varepsilon_i = 0$, $\mathbb{E}\varepsilon_i^2 = 1$ and $\beta_i = \mathbb{E}[\xi^4(i)]/[\mathbb{E}\xi^2(i)]^2 - 3 = -2$. This last value denotes a departure from the Gaussian case.

As in the previous example, we examine the limiting distributions of the three largest and the three smallest eigenvalues $\{\lambda_{n,j}, j = 1, \ldots, 6\}$ of the sample covariance matrix. Following Theorem 3.1, we have

- for $j = 1$ and 6,

$$\delta_{n,j} = \sqrt{n}[\lambda_{n,j} - \phi(\alpha_k)] \overset{\mathcal{D}}{\Longrightarrow} \mathcal{N}(0, s_{\alpha_k}^2), \quad s_{\alpha_k}^2 = \sigma_{\alpha_k}^2 \frac{y}{(\alpha_k - 1)^2}.$$

  Compared to the previous Gaussian case, as the factor $y/(\alpha_k - 1)^2 < 1$, the limiting Gaussian distributions of the largest and the smallest eigenvalue are less dispersed;

- for $j = (2, 3)$ or $j = (4, 5)$, the two-dimensional vector $\delta_{n,j} = \sqrt{n}[\lambda_{n,j} - \phi(\alpha_k)]$ converges weakly to the distribution of (ordered) eigenvalues of the random matrix

$$G = \sigma_{\alpha_k} \begin{pmatrix} W_{11} & W_{12} \\ W_{12} & W_{22} \end{pmatrix}. \tag{4.4}$$

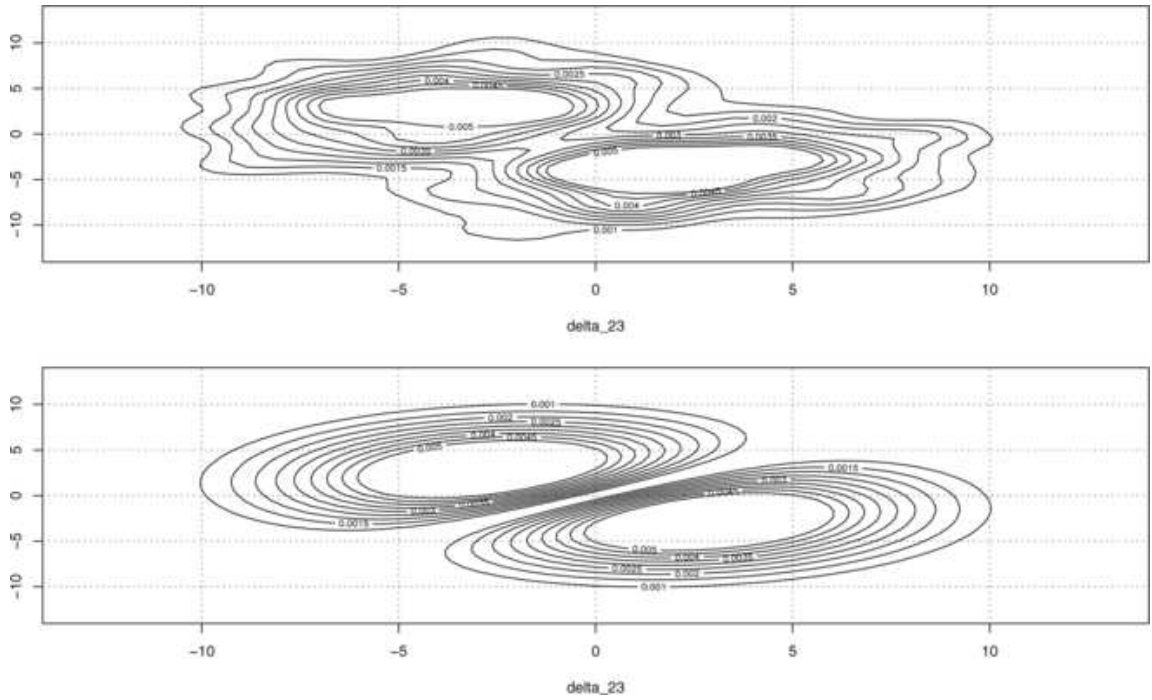

Fig. 2. Limiting bivariate distribution from the second and the third sample eigenvalues. Top: Contour lines of the empirical kernel density estimates $\widehat{f}_{n,23}$ from 1000 independent replications with $p = 500$, $n = 1000$ and Gaussian entries. Bottom: Contour lines of their limiting distribution given by the eigenvalues of a $2 \times 2$ Gaussian–Wigner matrix.

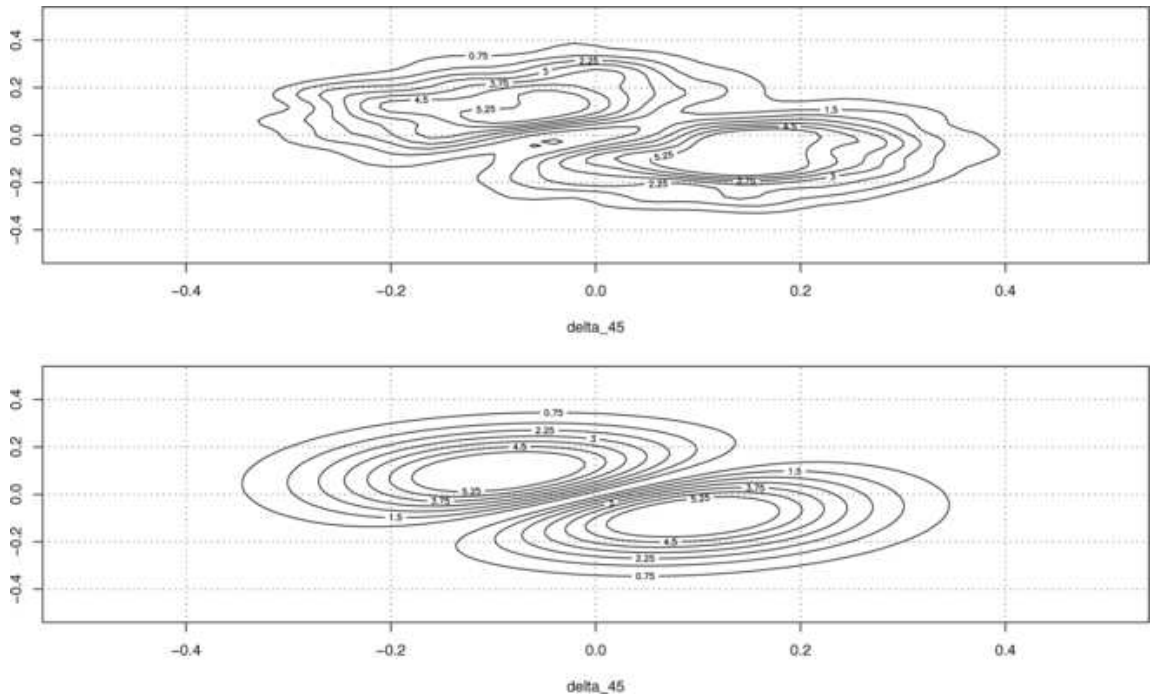

Fig. 3. Limiting bivariate distribution from the second and the third smallest sample eigenvalues. Top: Contour lines of the empirical kernel density estimates $\widehat{f}_{n,45}$ from 1000 independent replications with $p = 500$, $n = 1000$ and Gaussian entries. Bottom: Contour lines of their limiting distribution given by the eigenvalues of a $2 \times 2$ Gaussian–Wigner matrix.

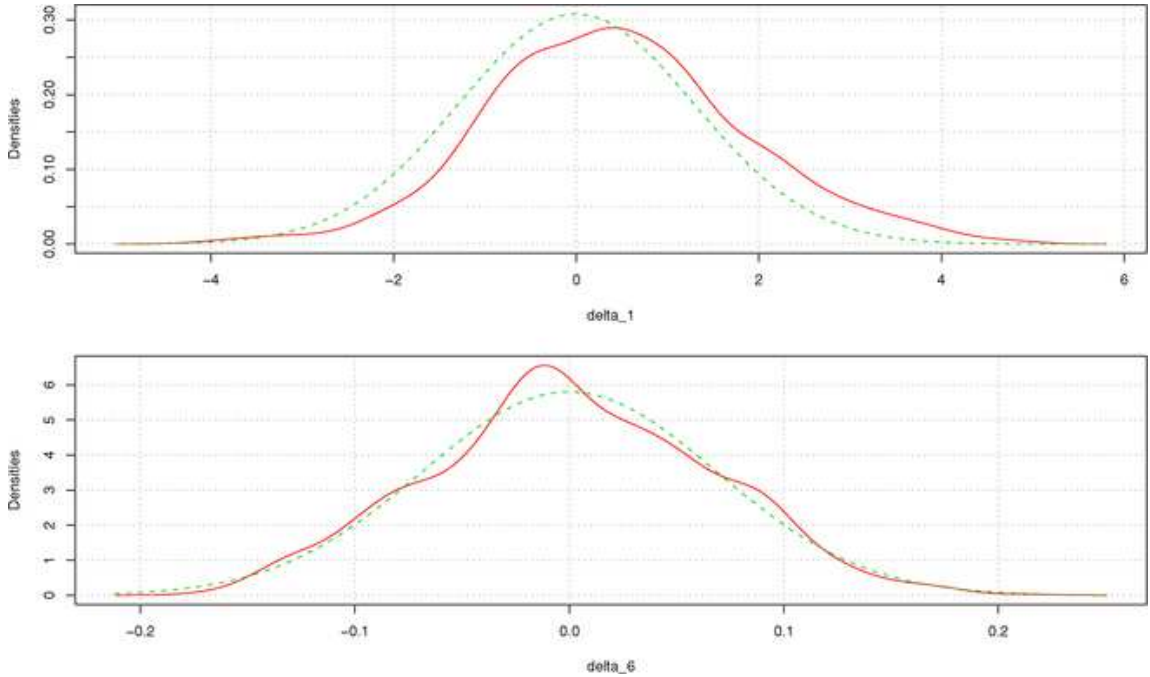

Fig. 4. Empirical density estimates (in solid lines) from the largest (top: $\widehat{f}_{n,1}$) and the smallest (bottom: $\widehat{f}_{n,6}$) sample eigenvalue from 1000 independent replications, compared to their Gaussian limits (dashed lines). *Binary entries* with $p = 500$ and $n = 1000$.

Here, because the initial variables $(\xi(i))$'s are binary, hence $\beta_i = -2$ (which is zero for Gaussian variables), by Eqs (3.11) and (3.12), we have $\mathrm{var}(W_{12}) = \frac{1}{2}$ but $\mathrm{var}(W_{11}) = \mathrm{var}(W_{22}) = y/(\alpha_k - 1)^2$. Therefore, the matrix $W = (W_{ij})$ is no more a real Gaussian–Wigner matrix. Again, the variance parameter $\sigma^2_{\alpha_k}$ is defined as previously but with $k = 2$ for $j = (2, 3)$ and $k = 3$ for $j = (4, 5)$, respectively. Unfortunately and unlike the previous Gaussian case, the joint distribution of eigenvalues of $W$ is unknown analytically. We then compute empirically by simulation this joint density using 10000 independent replications. Again, as $y/(\alpha_k - 1)^2 < 1$, these limiting distributions are less dispersed than previously.

The kernel density estimates $\widehat{f}_{n,1}$, $\widehat{f}_{n,6}$, $\widehat{f}_{n,23}$ and $\widehat{f}_{n,45}$ are computed as in the previous case using $p = 500$, $n = 1000$ and 1000 independent replications.

Figure 4 compares the two univariate density estimates $\widehat{f}_{n,1}$ and $\widehat{f}_{n,6}$ to their Gaussian limits. Again, we see that simulations confirm well the found formula. However, we remark a slower convergence rate than in the Gaussian case.

The bivariate density estimates $\widehat{f}_{n,23}$ and $\widehat{f}_{n,45}$ are then compared to their limiting densities in Figs 5 and 6, respectively. Again we see that the theoretical result is well confirmed. We remark that the shape of these bivariate limiting distributions is rather different from the previous Gaussian case. We remind the reader that the limiting bivariate densities are obtained by simulations of 10000 independent $G$ matrices given in (4.4).

## 5. Some extensions

It is possible to extend the spiked population model introduced in Section 2 to a much greater generality. Let us consider a population $p \times p$ covariance matrix

$$V = \mathrm{cov}(x) = \begin{pmatrix} \Sigma & 0 \\ 0 & T_p \end{pmatrix},$$

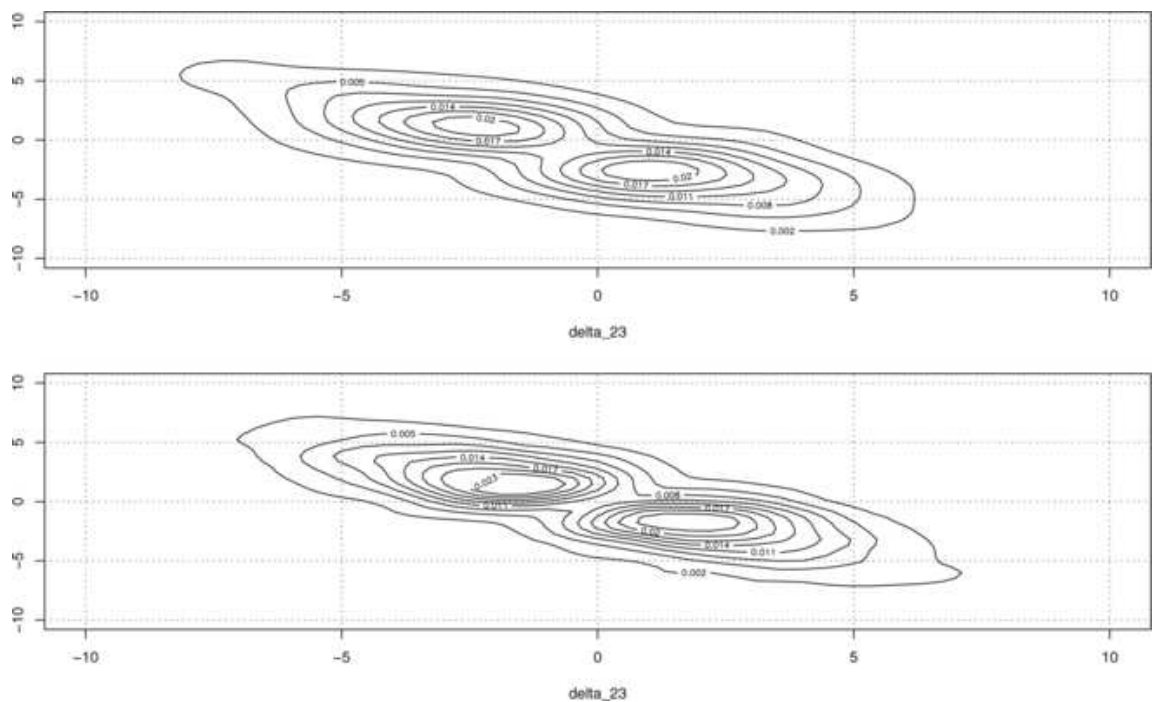

Fig. 5. Limiting bivariate distribution from the second and the third sample eigenvalues. Top: Contour lines of the empirical kernel density estimates $\widehat{f}_{n,23}$ from 1000 independent replications with $p = 500$, $n = 1000$ and *binary entries*. Bottom: Contour lines of their limiting distribution given by the eigenvalues of a $2 \times 2$ random matrix (computed by simulations).

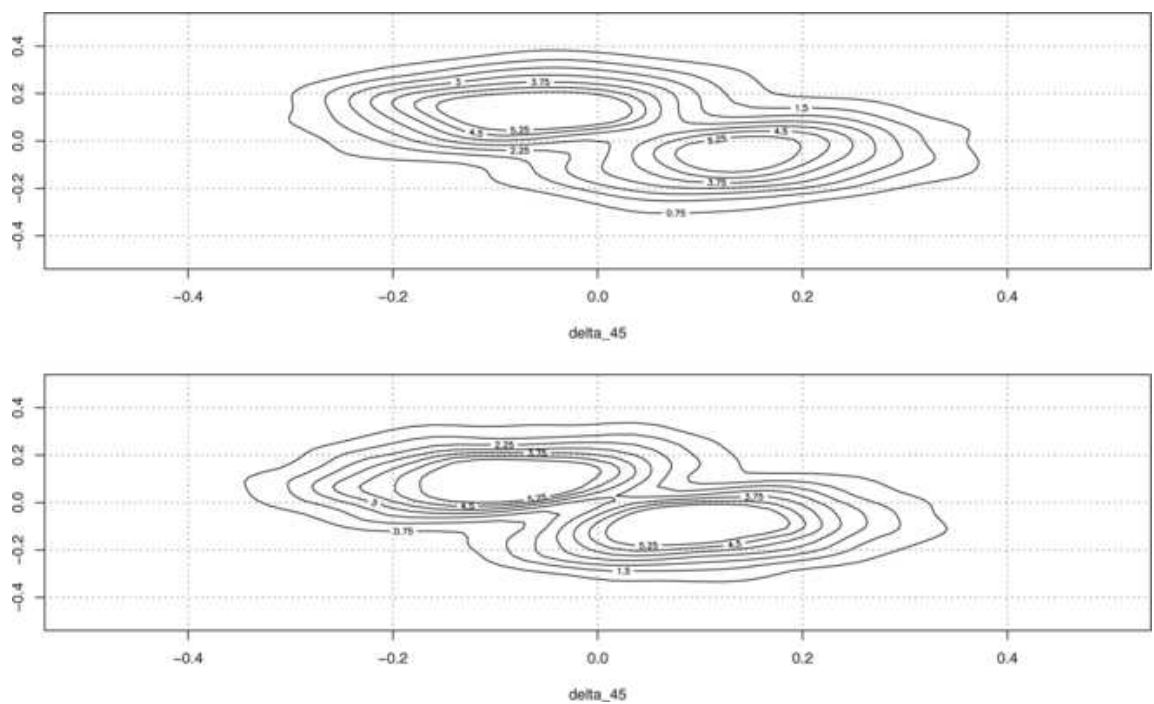

Fig. 6. Limiting bivariate distribution from the second and the third smallest sample eigenvalues. Top: Contour lines of the empirical kernel density estimates $\widehat{f}_{n,45}$ from 1000 independent replications with $p = 500$, $n = 1000$ and *binary entries*. Bottom: Contour lines of their limiting distribution given by the eigenvalues of a $2 \times 2$ random matrix (computed by simulations).

where $\Sigma$ is as previously while $T_p$ is now an arbitrary Hermitian matrix. As $M$ will be fixed and $p \to \infty$, the limit $F$ of the ESD of the sample covariance matrix depends on the sequence of $(T_p)$ only. With some ambiguity, we again call the eigenvalues $\alpha_k$'s of $\Sigma$ *spike eigenvalues* in the sense that they do not contribute to this limit.

In the following, we assume for simplicity that $\Sigma$ as well as $T_p$ are diagonal, and when $p \to \infty$, the empirical distribution of the eigenvalues of $T_p$ converges weakly to a probability measure $H(\mathrm{d}t)$ on the real line. Therefore, the limit $F$ of the ESD is characterized by an explicit formula for its Stieltjies transform, see [5].

The previous model of Section 2 corresponds to the situation where $H(\mathrm{d}t)$ is the Dirac measure at the point 1. A more involved example which will be analyzed later by numerical computations is the following. The core spectrum of $V$ is made with two eigenvalues $\omega_1 > \omega_2 > 0$, nearly $p/2$ times for each, and $V$ has a fixed number $M$ of spiked eigenvalues distinct from the $\omega_j$'s. In this case, the limiting distribution $H$ is $\frac{1}{2}(\delta_{\{\omega_1\}}(\mathrm{d}t) + \delta_{\{\omega_2\}}(\mathrm{d}t))$, a mixture of two Dirac masses.

The sample eigenvalues $\{\lambda_{n,j}\}$ are defined as previously. Assume that a spiked eigenvalue $\alpha_k$ is "sufficiently separated" from the core spectrum of $V$, so that for some function $\psi$ to be determined, there is a point $\psi(\alpha_k)$ outside the support of $F$ to which converge almost surely $n_k$ packed sample eigenvalues $\{\lambda_{n,j}, j \in J_k\}$. In such a case, the analysis we have proposed is also valid yielding a CLT analogous to Theorem 3.1: the $n_k$-dimensional real vector

$$\sqrt{n}\{\lambda_{n,j} - \psi(\alpha_k), j \in J_k\}$$

converges weakly to the distribution of the $n_k$ eigenvalues of some Gaussian random matrix. In particular, if $n_k = 1$, this limiting distribution is Gaussian.

We do not intend to provide here all details in this extended situation. However, let us indicate how we can determine the almost sure limit $\psi(\alpha_k)$ of the packed eigenvalues. From the almost sure convergence and since $\psi(\alpha_k)$ is outside the support of $F$, with probability tending to one, $\lambda_{n,j}$ solve the determinant equation (3.3). With $A_n = X_2^*(\lambda I - X_2 X_2^*)^{-1}X_2$, we have

$$K_n(\lambda) = S_{11} + X_1 A_n X_1^* = \frac{1}{n}\xi_{1:n}(I + A_n)\xi_{1:n}^*,$$

which tends almost surely to $[1 + ym_1(\lambda)]\Sigma$. Therefore, any limit $\lambda$ of a $\lambda_{n,j}$ fulfills the relation

$$\lambda - [1 + ym_1(\lambda)]\alpha = 0, \tag{5.1}$$

for some eigenvalue $\alpha$ of $\Sigma$.

Let $m(\lambda)$ be the Stieltjies transform of the limiting distribution $F$ and $\underline{m}(\lambda)$ the one of $yF(\mathrm{d}t) + (1-y)\delta_{\{0\}}(\mathrm{d}t)$. Clearly, $\lambda\underline{m}(\lambda) = -1 + y + y\lambda m(\lambda)$. Moreover, it is known that, see e.g. [5],

$$\lambda = -\frac{1}{\underline{m}(\lambda)} + y\int\frac{t}{1 + t\underline{m}(\lambda)}H(\mathrm{d}t). \tag{5.2}$$

As $m_1(\lambda) = -1 - \lambda m(\lambda)$ by definition, Eq. (5.1) reads as

$$\lambda = [1 - y - y\lambda m(\lambda)]\alpha = -\lambda\underline{m}(\lambda)\alpha.$$

It follows then $1 + \alpha\underline{m}(\lambda) = 0$ (generally, $\lambda \neq 0$). Combining with (5.2), we get finally

$$\lambda = \psi(\alpha) = \alpha\left[1 + y\int\frac{t}{\alpha - t}H(\mathrm{d}t)\right]. \tag{5.3}$$

In particular, for the original spiked population model with $H(\mathrm{d}t) = \delta_{\{1\}}(\mathrm{d}t)$, we recover the relation given in (2.3).

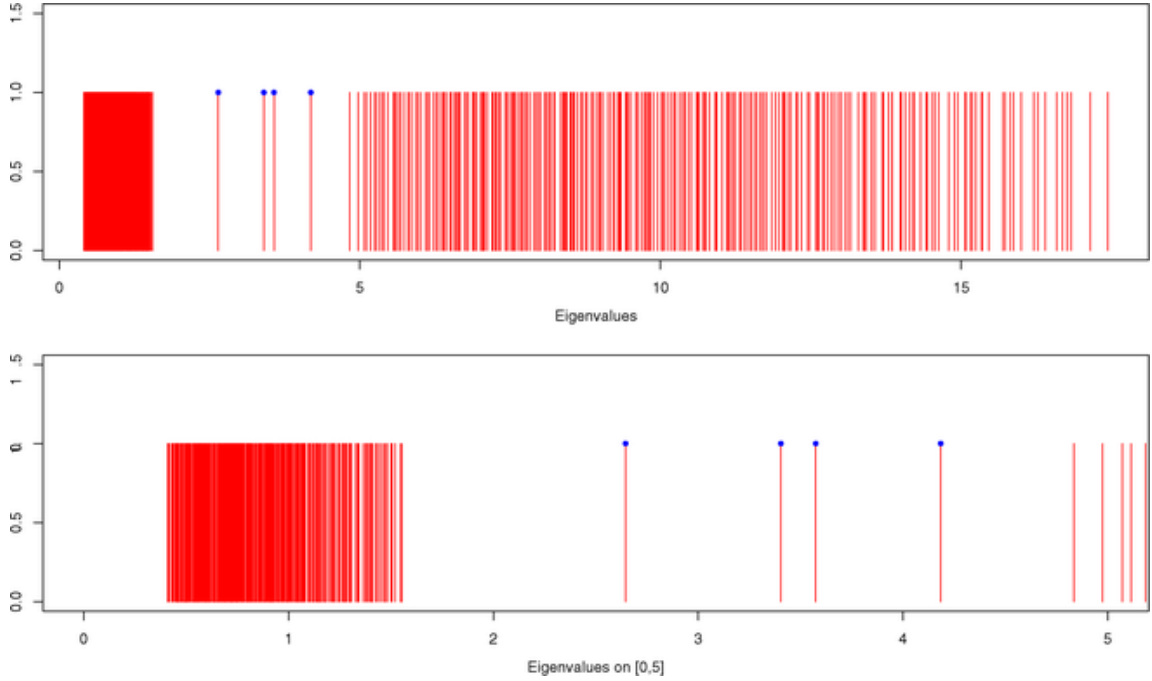

Fig. 7. An example of $p = 500$ sample eigenvalues (top) and a zoomed view on $[0,5]$ (bottom). The limiting distribution of the ESD has support $[0.395, 1.579] \cup [4.784, 17.441]$. The four eigenvalues $\{\lambda_{n,j}, 1 \leq \lambda \leq 4\}$ in the middle, related to spiked eigenvalues, are marked with a point. Gaussian entries with $n = 2500$.

We conclude the section by giving some numerical results of the above mentioned example of an extended spiked population model. Then, we consider $(\omega_1, \omega_2) = (1, 10)$, $(\alpha_1, \alpha_2, \alpha_3) = (5, 4, 3)$ with respective multiplicity $(1, 2, 1)$, and the limit ratio $y = 0.2$. Note that these spiked eigenvalues are now *between* the dominating eigenvalues (1 and 10). On the other hand, the support of the limiting distribution of the ESD can be determined following the method given in [5], and we get two disjointed intervals: $\operatorname{supp} F = [0.395, 1.579] \cup [4.784, 17.441]$.

For simulation, we use $p = 500$, $n = 2500$ and the eigenvalues of the population covariance matrix $V$ are 1 (248), 3 (1), 4 (2), 5 (1) and 10 (248). We simulate 500 independent replications of the sample covariance matrix with Gaussian variables. An example of these 500 replications is displayed in Fig. 7.

For each replication, the four eigenvalues at the middle (of indexes 249, 250, 251, 252) are extracted. Let us denote these 4 eigenvalues by $\lambda_{n,1}, \lambda_{n,2}, \lambda_{n,3}, \lambda_{n,4}$. By (5.3), we know that the almost sure limits of these sample eigenvalues are respectively

$$\psi(\alpha_k) = \alpha_k \left[ 1 + \frac{1}{10(\alpha_k - 1)} + \frac{1}{\alpha_k - 10} \right] = (4.125, 3.467, 2.721).$$

The next Fig. 8 displays the empirical densities of

$$\delta_{n,j} = \sqrt{n}(\lambda_{n,j} - \psi(\alpha_k)), \quad 1 \leq j \leq 4,$$

from the 500 independent replications. The graphs of $\delta_{n,1}$ and $\delta_{n,4}$ confirm a limiting zero-mean Gaussian distribution corresponding to single spike eigenvalues 5 and 3. On the contrary, the limiting distributions of $\delta_{n,2}$ and $\delta_{n,3}$, related to the double spike eigenvalue 4, are not zero-mean Gaussian. We note that $\delta_{n,2}$ and $-\delta_{n,3}$ have approximately the same distribution. Indeed, their joint distribution converges to that of the eigenvalues of $2 \times 2$ Gaussian–Wigner matrix.

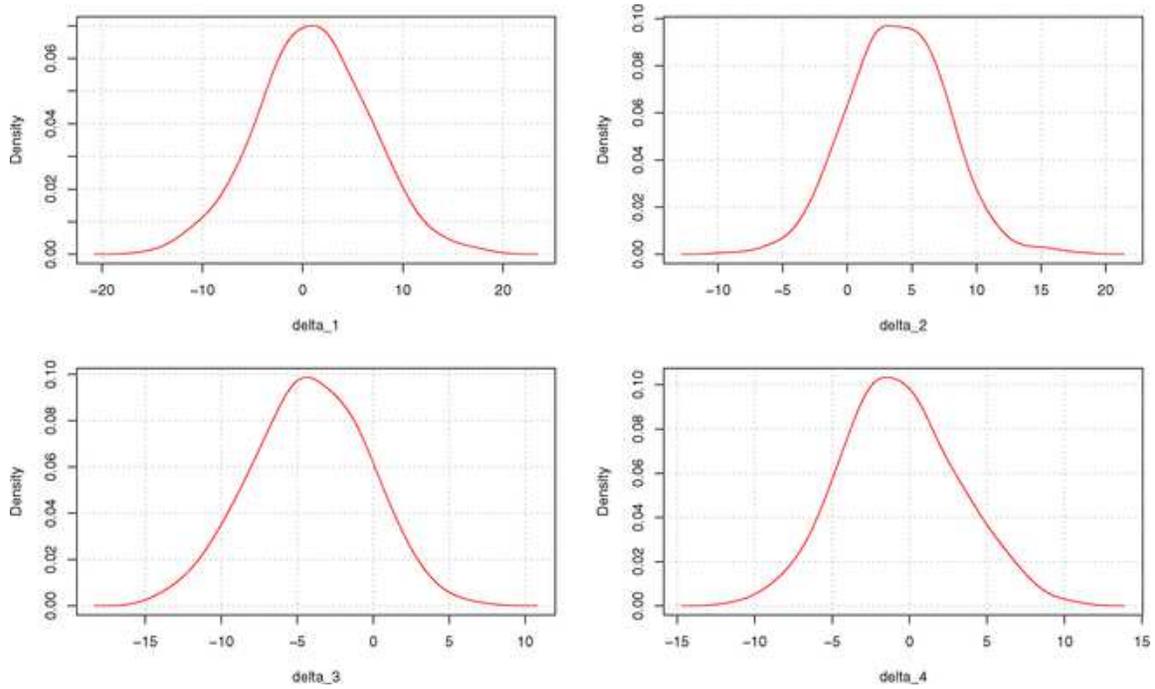

Fig. 8. Empirical densities of the normalized sample eigenvalues $\{\delta_{n,j}, 1 \leq j \leq 4\}$ from 500 independent replications. Gaussian entries with $p = 500$ and $n = 2500$.

## 6. Proofs of Propositions 3.1, 3.2 and Theorem 3.1

Before giving the proofs, some preliminary results and useful lemmas are introduced. Note that these proofs are based on a CLT for random sesquilinear forms which is itself introduced and proved in Section 7.

### 6.1. Preliminary results and useful lemmas

For $\lambda \notin [a_y, b_y]$, we define

$$m_1(\lambda) = \int \frac{x}{\lambda - x} F_y(\mathrm{d}x), \tag{6.1}$$

$$m_2(\lambda) = \int \frac{x^2}{(\lambda - x)^2} F_y(\mathrm{d}x), \tag{6.2}$$

$$m_3(\lambda) = \int \frac{x}{(\lambda - x)^2} F_y(\mathrm{d}x). \tag{6.3}$$

It is easily seen that

$$\int \frac{\lambda}{\lambda - x} F_y(\mathrm{d}x) = 1 + m_1(\lambda), \qquad \int \frac{\lambda^2}{(\lambda - x)^2} = 1 + 2m_1(\lambda) + m_2(\lambda).$$

If a real constant $\alpha \notin [1 - \sqrt{y}, 1 + \sqrt{y}]$, then $\phi(\alpha) \notin [a_y, b_y]$ and we have

$$m_1 \circ \phi(\alpha) = \frac{1}{\alpha - 1}, \tag{6.4}$$

$$m_2 \circ \phi(\alpha) = \frac{(\alpha - 1) + y(\alpha + 1)}{(\alpha - 1)[(\alpha - 1)^2 - y]}, \tag{6.5}$$

$$m_3 \circ \phi(\alpha) = \frac{1}{(\alpha-1)^2 - y}. \tag{6.6}$$

Let us mention that all these formulas can be obtained by derivation of the Stieltjes transform of the Marčenko–Pastur law $F_y(\mathrm{d}x)$

$$m(z) = \int \frac{1}{x-z} F_y(\mathrm{d}x) = \frac{1}{2yz}\{1 - y - z + \sqrt{(y+1-z)^2 - 4y}\}, \quad z \notin [a_y, b_y].$$

Here, $\sqrt{u}$ denotes the square root with positive imaginary part for $u \in \mathbb{C}$.

The following lemma gives the law of large numbers for some useful statistics related to the random matrix $A_n$ introduced in Eq. (3.4).

**Lemma 6.1.** *We have*

$$\frac{1}{n}\operatorname{tr} A_n \xrightarrow{P} ym_1(\lambda), \tag{6.7}$$

$$\frac{1}{n}\operatorname{tr} A_n A_n^* \xrightarrow{P} ym_2(\lambda), \tag{6.8}$$

$$\frac{1}{n}\sum_{i=1}^n a_{ii}^2 \xrightarrow{P} \left(\frac{y[1+m_1(\lambda)]}{\lambda - y[1+m_1(\lambda)]}\right)^2. \tag{6.9}$$

**Proof.** Let $\beta_{n,j}, j = 1, \ldots, p$ be the eigenvalues of $S_{22} = X_2 X_2^*$. The first equality is easy. For the second one, we have

$$\frac{1}{n}\operatorname{tr} A_n A_n^* = \frac{1}{n}\operatorname{tr}(\lambda I - X_2 X_2^*)^{-1} X_2 X_2^* (\lambda I - X_2 X_2^*)^{-1} X_2 X_2^*$$

$$= \frac{p}{n}\sum_{j=1}^p \frac{\beta_{n,j}^2}{(\lambda - \beta_{n,j})^2} \xrightarrow{P} y\int \frac{x^2}{(\lambda-x)^2} F_y(\mathrm{d}x).$$

For (6.9), let $e_i \in \mathbb{C}^n$ be the column vector whose $i$th element is 1 and others are 0 and $X_{2i}$ denote the matrix obtained from $X_2$ by deleting the $i$th column of $X_2$. We have $X_2 = X_{2i} + \frac{1}{n}\eta_i\eta_i^*$. Therefore,

$$a_{ii} = e_i^* X_2^*(\lambda I - X_2 X_2^*)^{-1} X_2 e_i = \frac{1}{n}\eta_i^*(\lambda I - X_2 X_2^*)^{-1}\eta_i = -\frac{\frac{1}{n}\eta_i^*(X_{2i}X_{2i}^* - \lambda I)^{-1}\eta_i}{1 + \frac{1}{n}\eta_i^*(X_{2i}X_{2i}^* - \lambda I)^{-1}\eta_i}.$$

Using Lemma 2.7 of [5],

$$\mathbb{E}\left|\frac{1}{n}\eta_i^*(X_{2i}X_{2i}^* - \lambda I)^{-1}\eta_i - \frac{1}{n}\operatorname{tr}(X_{2i}X_{2i}^* - \lambda I)^{-1}\right|^2 \le \frac{K}{n^2}\mathbb{E}|\eta(1)|^4 \mathbb{E}\operatorname{tr}(X_{2i}X_{2i}^* - \lambda I)^{-2},$$

which gives that

$$a_{ii} \xrightarrow{P} -\frac{y\int 1/(x-\lambda)F_y(\mathrm{d}x)}{1 + y\int 1/(x-\lambda)F_y(\mathrm{d}x)} = \frac{y[1+m_1(\lambda)]}{\lambda - y[1+m_1(\lambda)]}. \tag{6.10}$$

Further, it is easy to verify that

$$\lim_{n\to\infty} \mathbb{E}\frac{\operatorname{tr}(X_2^*(\lambda I - X_2 X_2^*)^{-1} X_2)^4}{n} < \infty,$$

which implies, together with inequality 3.3.41 of [8] that

$$\sup_n \mathbb{E}a_{11}^4 = \sup_n \frac{1}{n}\sum_{i=1}^n \mathbb{E}a_{ii}^4 \le \sup_n \mathbb{E}\frac{\operatorname{tr}(X_2^*(\lambda I - X_2 X_2^*)^{-1} X_2)^4}{n} < \infty.$$

Therefore, the family of the random variables $\{a_{11}^2\}$ indexed by $n$ is uniformly integrable. Combining with (6.10), we get

$$\mathbb{E}\left|\frac{1}{n}\sum_{i=1}^{n}a_{ii}^2 - \left(\frac{y[1+m_1(\lambda)]}{\lambda - y[1+m_1(\lambda)]}\right)^2\right| \leq \mathbb{E}\left|a_{11}^2 - \left(\frac{y[1+m_1(\lambda)]}{\lambda - y[1+m_1(\lambda)]}\right)^2\right| \to 0.$$

Thus (6.9) follows. □

### 6.2. Proof of Proposition 3.1

We apply Theorem 7.1 by considering $K = \frac{1}{2}M(M+1)$ bilinear forms

$$u(i)(I+A_n)u(j)^{\mathrm{T}}, \quad 1 \leq i \leq j \leq M,$$

with

$$u(i) = (\xi_1(i), \ldots, \xi_n(i)).$$

More precisely, with $\ell = (i,j)$, we are substituting $u(i)^{\mathrm{T}}$ for $X(\ell)$, and $u(j)^{\mathrm{T}}$ for $Y(\ell)$, respectively. Consequently, $x_{\ell 1} = \xi_1(i)$ and $y_{\ell 1} = \xi_1(j)$ for the application of Theorem 7.1.

We have, by Lemma 6.1,

$$\theta = \tau = \lim_{n}\frac{1}{n}\operatorname{tr}(I+A_n)^2 = 1 + 2ym_1(\lambda) + ym_2(\lambda),$$

$$\omega = \lim_{n}\frac{1}{n}\sum_{i=1}^{n}[(I+A_n)_{ii}]^2 = 1 + 2ym_1(\lambda) + \left(\frac{y[1+m_1(\lambda)]}{\lambda - y[1+m_1(\lambda)]}\right)^2.$$

Following Theorem 7.1, $R_n$ converges weakly to a symmetric random matrix with zero-mean Gaussian variables $R = (R_{ij})$ with the following covariance function, assuming $1 \leq i \leq j \leq M$,

$$\begin{aligned}
\operatorname{cov}(R_{ij}, R_{i'j'}) &= \omega\{\mathbb{E}[\xi(i)\xi(j)\xi(i')\xi(j')] - \Sigma_{ij}\Sigma_{i'j'}\} + (\theta - \omega)\{\mathbb{E}[\xi(i)\xi(j')]\mathbb{E}[\xi(i')\xi(j)]\} \\
&\quad + (\theta - \omega)\{\mathbb{E}[\xi(i)\xi(i')]\mathbb{E}[\xi(j)\xi(j')]\}.
\end{aligned} \tag{6.11}$$

### 6.3. Proof of Proposition 3.2

The aim is to apply Theorem 7.3 to $K = \frac{1}{2}M(M+1)$ sesquilinear forms

$$u(i)(I+A_n)u(j)^*, \quad 1 \leq i \leq j \leq M,$$

with

$$u(i) = (\xi_1(i), \ldots, \xi_n(i)).$$

More precisely, with $\ell = (i,j)$, we are substituting $u(i)^*$ for $X(\ell)$, and $u(j)^*$ for $Y(\ell)$, respectively. Consequently, $x_{\ell 1} = \bar{\xi}_1(i)$ and $y_{\ell 1} = \bar{\xi}_1(j)$ for the application of Theorem 7.3.

Again by Lemma 6.1,

$$\theta = \lim_{n}\frac{1}{n}\operatorname{tr}(I+A_n)^2 = 1 + 2ym_1(\lambda) + ym_2(\lambda),$$

$$\omega = \lim_{n}\frac{1}{n}\sum_{i=1}^{n}[(I+A_n)_{ii}]^2 = 1 + 2ym_1(\lambda) + \left(\frac{y[1+m_1(\lambda)]}{\lambda - y[1+m_1(\lambda)]}\right)^2.$$

Here we need an additional condition which is specific to the complex case. Assume therefore,

$$\tau = \lim_n \frac{1}{n} \operatorname{tr}(I + A_n)(I + A_n)^{\mathrm{T}} = 1 + 2ym_1(\lambda) + m_4(\lambda).$$

Consequently by Theorem 7.3, $R_n$ converges weakly to a zero-mean Hermitian random matrix $R = (R_{ij})$. Moreover, the joint distribution of the real and imaginary parts of the upper-triangular bloc $\{R_{ij}, 1 \leq i \leq j \leq M\}$ is a $2K$-dimensional Gaussian vector with covariance matrix

$$\Gamma = \begin{pmatrix} \Gamma_{11} & \Gamma_{12} \\ \Gamma_{21} & \Gamma_{22} \end{pmatrix}, \tag{6.12}$$

where

$$\Gamma_{11} = \frac{1}{4} \sum_{j=1}^{3} \{ 2\Re(B_j) + B_{ja} + B_{jb} \},$$

$$\Gamma_{22} = \frac{1}{4} \sum_{j=1}^{3} \{ -2\Re(B_j) + B_{ja} + B_{jb} \},$$

$$\Gamma_{12} = \frac{1}{2} \sum_{j=1}^{3} \Im(B_j),$$

and for $1 \leq i \leq j \leq M$ and $1 \leq i' \leq j' \leq M$, with the $B$-matrices defined in the proposition.

## 6.4. Proof of Theorem 3.1

Let $\alpha_k \notin [1 - \sqrt{y}, 1 + \sqrt{y}]$ be fixed. Following Section 3.1, we can assume that the $n_k$ packed sample eigenvalues $\{\lambda_{n,j}, j \in J_k\}$ are solutions of the equation $|\lambda - K_n(\lambda)| = 0$. As $\lambda_{n,j} \to \lambda_k$ almost surely, we define

$$\delta_{n,j} = \sqrt{n}(\lambda_{n,j} - \lambda_k).$$

We have

$$\lambda_{nj}I - K_n(\lambda_{n,j}) = \lambda_k I + \frac{1}{\sqrt{n}}\delta_{n,j}I - K_n(\lambda_k) - [K_n(\lambda_{n,j}) - K_n(\lambda_k)].$$

Furthermore, using $A^{-1} - B^{-1} = A^{-1}(B - A)B^{-1}$, we have

$$K_n(\lambda_{n,j}) - K_n(\lambda_k) = \frac{1}{n}\xi_{1:n}X_2^* \left\{ \left( \left[ \lambda_k + \frac{1}{\sqrt{n}}\delta_{n,j} \right] I - S_{22} \right)^{-1} - (\lambda_k I - S_{22})^{-1} \right\} X_2 \xi_{1:n}^*$$

$$= -\frac{1}{\sqrt{n}}\delta_{n,j}\frac{1}{n}\xi_{1:n}X_2^* \left( \left[ \lambda_k + \frac{1}{\sqrt{n}}\delta_{n,j} \right] I - S_{22} \right)^{-1} (\lambda_k I - S_{22})^{-1} X_2 \xi_{1:n}^*$$

$$= -\frac{1}{\sqrt{n}}\delta_{n,j}[ym_3(\lambda_k)\Sigma + \mathrm{o}_P(1)].$$

Combining these estimations and (3.5), (3.6), we have

$$\lambda_{nj}I - K_n(\lambda_{n,j}) = \lambda_k I - [1 + ym_1(\lambda_k)]\Sigma - \frac{1}{\sqrt{n}}R_n(\lambda_k) + \frac{1}{\sqrt{n}}\delta_{n,j}[I + ym_3(\lambda_k)\Sigma] + \mathrm{o}_P\left( \frac{1}{\sqrt{n}} \right). \tag{6.13}$$

By Section 3.1, $R_n(\lambda_k)$ converges in distribution to a $M \times M$ random matrix $R(\lambda_k)$ with Gaussian entries with a fully identified covariance matrix. We now follow a method devised in [1] and [2] for limiting distributions of eigenvalues or eigenvectors from random matrices. First, we use Skorokhod strong representation so that on an appropriate probability space, the convergence $R_n(\lambda_k) \to R(\lambda_k)$ as well as (6.13) take place almost surely. Multiplying both sides of (6.13) by $U$ from the left and by $U^*$ from the right yields

$$U[\lambda_{nj}I - K_n(\lambda_{n,j})]U^* = \begin{pmatrix} \ddots & 0 & 0 \\ 0 & (\lambda_k - [1 + ym_1(\lambda_k)]\alpha_u)I_{n_u} & 0 \\ 0 & 0 & \ddots \end{pmatrix} - \frac{1}{\sqrt{n}}UR_n(\lambda_k)U^*$$

$$+ \frac{1}{\sqrt{n}}\begin{pmatrix} \ddots & 0 & 0 \\ 0 & \delta_{n,j}(1 + ym_3(\lambda_k)\alpha_u)I_{n_u} & 0 \\ 0 & 0 & \ddots \end{pmatrix} + o\left(\frac{1}{\sqrt{n}}\right).$$

First, in the right-hand side of the equation and using a bloc decomposition induced by (3.17), we see that all the non diagonal blocs tend to zero. Next, for a diagonal bloc with index $u \neq k$, by definition $\lambda_k - [1 + ym_1(\lambda_k)]\alpha_u \neq 0$, and this is the limit of that diagonal bloc since the contributions from the remaining three terms tend to zero. As $\lambda_k - [1 + ym_1(\lambda_k)]\alpha_k = 0$ by definition, the $k$th diagonal bloc reduces to

$$-\frac{1}{\sqrt{n}}[UR_n(\lambda_k)U^*]_{kk} + \frac{1}{\sqrt{n}}\delta_{n,j}(1 + ym_3(\lambda_k)\alpha_k)I_{n_k} + o\left(\frac{1}{\sqrt{n}}\right).$$

For $n$ sufficiently large, its determinant must be equal to zero,

$$\left| -\frac{1}{\sqrt{n}}[UR_n(\lambda_k)U^*]_{kk} + \frac{1}{\sqrt{n}}\delta_{n,j}(1 + ym_3(\lambda_k)\alpha_k)I_{n_k} + o\left(\frac{1}{\sqrt{n}}\right) \right| = 0,$$

or equivalently,

$$|-[UR_n(\lambda_k)U^*]_{kk} + \delta_{n,j}(1 + ym_3(\lambda_k)\alpha_k)I_{n_k} + o(1)| = 0.$$

Therefore, $\delta_{n,j}$ tends to a solution of

$$|-[UR_n(\lambda_k)U^*]_{kk} + \lambda(1 + ym_3(\lambda_k)\alpha_k)I_{n_k}| = 0,$$

that is, an eigenvalue of the matrix $(1 + ym_3(\lambda_k)\alpha_k)^{-1}\widetilde{R}_{kk}(\lambda_k)$. Finally, as the index $j$ is arbitrary, all the $J_k$ random variables $\sqrt{n}\{\lambda_{n,j} - \lambda_k, j \leq J_k\}$ converge almost surely to the set of eigenvalues of the above matrix. Of course, this convergence also holds in distribution on the new probability space, hence on the original one.

## 7. A CLT for random sesquilinear forms

The aim of this section is to establish a CLT for random sesquilinear forms as one of the central tools used in the paper. These results are independent from the previous sections and should have their own interest.

Consider a sequence $\{(x_i, y_i)_{i \in N}\}$ of i.i.d. complex-valued, zero-mean random vectors belonging to $\mathbb{C}^K \times \mathbb{C}^K$ with a finite moment of the fourth-order. We write

$$x_i = (x_{\ell i}) = \begin{pmatrix} x_{1i} \\ \vdots \\ x_{Ki} \end{pmatrix}, \qquad X(\ell) = (x_{\ell 1}, \ldots, x_{\ell n})^{\mathrm{T}}, \quad 1 \leq \ell \leq K, \tag{7.1}$$

with a similar definition for the vectors $\{Y(\ell)\}_{1 \leq \ell \leq K}$. Set $\rho(\ell) = \mathbb{E}[\bar{x}_{\ell 1}y_{\ell 1}]$.

**Theorem 7.1.** *Let $\{A_n = [a_{ij}(n)]\}_n$ be a sequence of $n \times n$ Hermitian matrices and the vectors $\{X(\ell), Y(\ell)\}_{1 \le \ell \le K}$ are as defined in (7.1). Assume that the following limits exist*

$$\omega = \lim_{n \to \infty} \frac{1}{n} \sum_{u=1}^{n} a_{uu}^2(n),$$

$$\theta = \lim_{n \to \infty} \frac{1}{n} \operatorname{tr} A_n^2 = \lim_{n \to \infty} \frac{1}{n} \sum_{u,v=1}^{n} |a_{uv}(n)|^2,$$

$$\tau = \lim_{n \to \infty} \frac{1}{n} \operatorname{tr} A_n A_n^{\mathrm{T}} = \lim_{n \to \infty} \frac{1}{n} \sum_{u,v=1}^{n} a_{uv}^2(n).$$

*Then, the $M$-dimensional complex-valued random vectors*

$$Z_n = (Z_{n,\ell}), \qquad Z_{n,\ell} = \frac{1}{\sqrt{n}} [X(\ell)^* A_n Y(\ell) - \rho(\ell) \operatorname{tr} A_n], \quad 1 \le \ell \le K, \tag{7.2}$$

*converge weakly to a zero-mean complex-valued vector $W$ whose real and imaginary parts are Gaussian. Moreover, the Laplace transform of $W$ is given by*

$$\mathbb{E} e^{c^{\mathrm{T}} W} = \exp\left[\frac{1}{2} c^{\mathrm{T}} B c\right], \quad c \in \mathbb{C}^K, \tag{7.3}$$

*where the matrix $B = B_1 + B_2 + B_3$ with*

$$
\begin{aligned}
B_1 &= \omega(\mathbb{E}[\bar{x}_{\ell 1} y_{\ell 1} \bar{x}_{\ell' 1} y_{\ell' 1}] - \rho(\ell)\rho(\ell')), & 1 \le \ell, \ell' \le K, \\
B_2 &= (\theta - \omega)(\mathbb{E}[\bar{x}_{\ell 1} y_{\ell' 1}] \mathbb{E}[\bar{x}_{\ell' 1} y_{\ell 1}]), & 1 \le \ell, \ell' \le K, \\
B_3 &= (\tau - \omega)(\mathbb{E}[\bar{x}_{\ell 1} \bar{x}_{\ell' 1}] \mathbb{E}[y_{\ell 1} y_{\ell' 1}]), & 1 \le \ell, \ell' \le K.
\end{aligned}
\tag{7.4}
$$

The proof of the theorem is postponed to the end of the section. First, we describe some specific applications of the theorem with their own interest. Note that by definition, the three matrices $B_j$'s are symmetrical (complex-valued).

Consider first the real variables case with i.i.d. random vectors $\{(x_i, y_i)_{i \in N}\}$ from $\mathbb{R}^K \times \mathbb{R}^K$, and a sequence of symmetric matrices $\{A_n = [a_{ij}(n)]\}_n$. We are then considering $K$ random bilinear forms and consequently, $\theta = \tau$. The matrix $B$ given above is then exactly the limiting covariance matrix of the Gaussian vector $W$.

**Corollary 7.1.** *Under the same conditions as in Theorem 7.1 but with real random vectors $\{(x_i, y_i)_{i \in N}\}$ and symmetric matrices $\{A_n\}_n$, the sequence of vectors $(Z_n)_n$ converges weakly to a zero-mean $K$-dimensional Gaussian vector with covariance matrix $B$.*

An interesting application to the case $(x_i) = (y_i)$ gives the following CLT for random quadratic forms in a straightforward way.

**Theorem 7.2.** *Let $\{A_n = [a_{ij}(n)]\}_n$ be a sequence of $n \times n$ real symmetric matrices, $(x_i)_{i \in \mathbb{N}}$ a sequence of i.i.d. $K$-dimensional real random vectors, with $\mathbb{E}[x_i] = 0$, $\mathbb{E}[x_1 x_1^{\mathrm{T}}] = (\gamma_{ij})$, $1 \le i, j \le K$, and $\mathbb{E}[\|x_1\|^4] < \infty$. Let the vectors $\{X(\ell)\}_{1 \le \ell \le K}$ be as defined in (7.1). Assume the following limits exist*

$$\omega = \lim_{n \to \infty} \frac{1}{n} \sum_{u=1}^{n} a_{uu}^2(n),$$

$$\theta = \lim_{n \to \infty} \frac{1}{n} \operatorname{tr} A_n^2.$$

*Then, the M-dimensional random vectors*

$$Z_n = (Z_{n,\ell}), \quad Z_{n,\ell} = \frac{1}{\sqrt{n}}[X(\ell)^{\mathrm{T}} A_n X(\ell) - \gamma_{\ell\ell} \operatorname{tr} A_n], \quad 1 \le \ell \le K, \tag{7.5}$$

*converge weakly to a zero-mean Gaussian vector with covariance matrix $D = D_1 + D_2$ where*

$$\begin{aligned}
D_1 &= \omega(\mathbb{E}[x_{\ell 1}^2 x_{\ell' 1}^2] - \gamma_{\ell\ell}\gamma_{\ell'\ell'}), \quad 1 \le \ell, \ell' \le K, \\
D_2 &= (\theta - \omega)(\gamma_{\ell\ell'}\gamma_{\ell'\ell} + \gamma_{\ell\ell'}^2), \quad 1 \le \ell, \ell' \le K.
\end{aligned} \tag{7.6}$$

If all the diagonal elements of the matrices $(A_n)$ are null, then $\omega = 0$. The limiting covariance matrix $D$ takes a much simpler form:

$$D = \theta(\gamma_{\ell\ell'}\gamma_{\ell'\ell} + \gamma_{\ell\ell'}^2), \quad 1 \le \ell, \ell' \le K.$$

For the general complex case, we need a special device. Write $Z_n = U_n + \mathrm{i}V_n$. Following Theorem 7.1, $(U_n, V_n)$ converges weakly to a $2K$-dimensional Gaussian vector with some covariance matrix $\Gamma$. The aim is to identify $\Gamma$. We have

$$\mathbb{E}\exp[t^{\mathrm{T}}U_n + s^{\mathrm{T}}V_n] \to \exp\left[\frac{1}{2}(t^{\mathrm{T}}, s^{\mathrm{T}})\Gamma\begin{pmatrix} t \\ s \end{pmatrix}\right], \quad t \in \mathbb{R}^K, s \in \mathbb{R}^K. \tag{7.7}$$

On the other hand, from $U_n = \frac{1}{2}(Z_n + \bar{Z}_n)$ and $V_n = \frac{1}{2\mathrm{i}}(Z_n - \bar{Z}_n)$, we have a second expression

$$\mathbb{E}\exp[t^{\mathrm{T}}U_n + s^{\mathrm{T}}V_n] = \mathbb{E}\exp\left[\left(\frac{t}{2} + \frac{s}{2\mathrm{i}}\right)^{\mathrm{T}}Z_n + \left(\frac{t}{2} - \frac{s}{2\mathrm{i}}\right)^{\mathrm{T}}\bar{Z}_n\right].$$

Interestingly enough, the last transform can be found by application of Theorem 7.1 to the random sesquilinear forms

$$\widetilde{Z}_n = \begin{pmatrix} Z_n \\ \overline{Z}_n \end{pmatrix}. \tag{7.8}$$

For ease of the presentation, we need to define more limiting quantities. For $1 \le \ell \le K$, let $\sigma_{X,\ell}^2 = \mathbb{E}[|x_{\ell 1}|^2], \sigma_{Y,\ell}^2 = \mathbb{E}[|y_{\ell 1}|^2]$. We introduce the following matrices

$$\begin{aligned}
B_{1a} &= \omega(\mathbb{E}[|x_{\ell 1}|^2|x_{\ell' 1}|^2] - \sigma_{X,\ell}^2\sigma_{X,\ell'}^2), \\
B_{1b} &= \omega(\mathbb{E}[|y_{\ell 1}|^2|y_{\ell' 1}|^2] - \sigma_{Y,\ell}^2\sigma_{Y,\ell'}^2), \\
B_{2a} &= (\theta - \omega)(|\mathbb{E}[\bar{x}_{\ell 1}x_{\ell' 1}]|^2), \\
B_{2b} &= (\theta - \omega)(|\mathbb{E}[\bar{y}_{\ell 1}y_{\ell' 1}]|^2), \\
B_{3a} &= (\tau - \omega)(|\mathbb{E}[x_{\ell 1}x_{\ell' 1}]|^2), \\
B_{3b} &= (\tau - \omega)(|\mathbb{E}[y_{\ell 1}y_{\ell' 1}]|^2).
\end{aligned} \tag{7.9}$$

Here, the indices are $1 \le \ell, \ell' \le K$. By definition, all these matrices are real and symmetrical. Let us also define the $2K \times 2K$ matrices

$$\widetilde{B}_j = \begin{pmatrix} B_j & B_{ja} \\ B_{jb} & \overline{B}_j \end{pmatrix}, \quad j = 1, 2, 3. \tag{7.10}$$

**Theorem 7.3.** *Consider the $M$-dimensional complex-valued random vectors $Z_n = (Z_{n,\ell})$ defined in Theorem 7.1. Under the the same conditions as in that theorem, the real and the imaginary parts $(U_n, V_n)$ of $Z_n$ converge weakly to a $2K$-dimensional Gaussian vector with covariance matrix*

$$\Gamma = \begin{pmatrix} \Gamma_{11} & \Gamma_{12} \\ \Gamma_{21} & \Gamma_{22} \end{pmatrix}, \tag{7.11}$$

*with*

$$\Gamma_{11} = \frac{1}{4} \sum_{j=1}^{3} \{2\Re(B_j) + B_{ja} + B_{jb}\},$$

$$\Gamma_{22} = \frac{1}{4} \sum_{j=1}^{3} \{-2\Re(B_j) + B_{ja} + B_{jb}\},$$

$$\Gamma_{12} = \frac{1}{2} \sum_{j=1}^{3} \Im(B_j).$$

**Proof.** For the vector of sesquilinear forms $\widetilde{Z}_n$ in (7.8), one can check that the limiting matrix $B$ in Theorem 7.1 is to be replaced by

$$\widetilde{B} = \begin{pmatrix} \sum_{j=1}^{3} B_j & \sum_{j=1}^{3} B_{ja} \\ \sum_{j=1}^{3} B_{jb} & \sum_{j=1}^{3} \overline{B}_j. \end{pmatrix}.$$

Then following this theorem, for $\widetilde{c} = (\frac{t^{\mathrm{T}}}{2} + \frac{s^{\mathrm{T}}}{2\mathrm{i}}, \frac{t^{\mathrm{T}}}{2} - \frac{s^{\mathrm{T}}}{2\mathrm{i}})^{\mathrm{T}}$,

$$\mathbb{E} \exp[\widetilde{c}^{\mathrm{T}} \widetilde{Z}_n] \to \exp\left\{\frac{1}{2}\left(\frac{t^{\mathrm{T}}}{2} + \frac{s^{\mathrm{T}}}{2\mathrm{i}}, \frac{t^{\mathrm{T}}}{2} - \frac{s^{\mathrm{T}}}{2\mathrm{i}}\right)\widetilde{B}\begin{pmatrix} t/2 + s/(2\mathrm{i}) \\ t/2 - s/(2\mathrm{i}) \end{pmatrix}\right\}.$$

By identifying this formula to Eq. (7.7), we get the required form of $\Gamma$. $\qquad \square$

### 7.1. *Proof of Theorem 7.1*

It is sufficient to establish the CLT for the sequence of linear combinations of random Hermitian forms

$$\sum_{\ell=1}^{K} c_\ell X(\ell)^* A_n Y(\ell),$$

where the coefficients $(c_\ell) \in \mathbb{C}^K$ are arbitrary. Notice that $\mathbb{E}[X(\ell)^* A_n Y(\ell)] = \rho(\ell)\operatorname{tr}(A_n)$, where $\rho(\ell) = \mathbb{E}[\bar{x}_{\ell 1} y_{\ell 1}]$.

First, by a classical procedure of truncation and renormalization (see Section 7.2 for details), we can, without loss of generality, assume that there is a sequence $\varepsilon_n \downarrow 0$ such that

$$1 \leq i \leq n, \quad \|x_i\| \vee \|y_i\| \leq \varepsilon_n n^{1/4}. \tag{7.12}$$

We will use the method of moments. Define, while dropping the index $n$ in the coefficients of $A_n$,

$$\xi_n = \frac{1}{\sqrt{n}} \sum_{\ell=1}^{K} c_\ell [X(\ell)^* A_n Y(\ell) - \rho(\ell)\operatorname{tr} A_n] = \frac{1}{\sqrt{n}} \sum_e a_e \psi_e,$$

where $e = (u, v) \in \{1, \ldots, n\}^2$ and

$$\psi_e = \begin{cases} \sum_{\ell=1}^{K} c_\ell [\bar{x}_u y_u - \rho(\ell)], & e = (u, u), \\ \sum_{\ell=1}^{K} c_\ell \bar{x}_u y_v, & e = (u, v), \ u \neq v. \end{cases}$$

Let $k \geq 1$ be a given integer. We have

$$n^{k/2} \xi_n^k = \sum_{e_1, \ldots, e_k} a_{e_1} \cdots a_{e_k} \psi_{e_1} \cdots \psi_{e_k}.$$

To each term in the sum we associate a directed graph $G$ by drawing an arrow $u \to v$ for each factor $e_j = (u, v)$. The set of vertices is then a subset of $\{1, \ldots, n\}$. Therefore, to a loop $u \to u$ corresponds the product $a_{uu} \psi_{u \to u}(\ell) = a_{uu} \sum_{\ell=1}^{K} c_\ell [\bar{x}_u(\ell) y_u(\ell) - \rho(\ell)]$ and to an edge $u \to v$ with $u \neq v$ corresponds the product $a_{uv} \psi_{u \to v} = a_{uv} \sum_{\ell=1}^{K} c_\ell \bar{x}_u(\ell) y_v(\ell)$. In other words,

$$n^{k/2} \xi_n^k = \sum_G a_G \psi_G, \quad a_G = \prod_{e \in G} a_e, \psi_G = \prod_{e \in G} \psi_e.$$

We now consider the collection of connected sub-graphs of $G$. These connected sub-graphs can be classified into two types.

- *Type-I sub-graphs.* We call $C$ a Type-I connected sub-graph of $G$ if $C$ contains loops only. In particular $C$ has a unique vertex. The set of all the $m_1$ Type-I connected sub-graphs is denoted by $\mathcal{F}_1$, and the degrees of their vertexes by $\mu_1, \ldots, \mu_{m_1}$, respectively.

  If $\mu_j = 2$ for some vertex $j$ in a sub-graph $C$, then $\mathbb{E} a_G \psi_G = 0$ because of independence. Therefore we need only consider those graphs $G$ whose $m_1$ Type-I sub-graphs have all their vertices with degrees $\mu_j \geq 4$. The contributions from all these sub-graphs to the moment part $\psi_G$ are then bounded by

$$\left| \mathbb{E} \prod_{C \in \mathcal{F}_1} \psi_C \right| \leq K(\varepsilon_n n^{1/4})^{\sum_{i=1}^{m_1} (\mu_i - 4)}. \tag{7.13}$$

- *Type-II sub-graphs.* A connected sub-graph containing at least one arrow $u \to v$ with $u \neq v$ is called a Type-II sub-graph. The set of all these $m_2$ components is denoted by $\mathcal{F}_2$. For each $C_s \in \mathcal{F}_2$, let $u_s$ be the number of its vertices whose degrees are denoted by $\gamma_{js}$, $j = 1, \ldots, u_s$. As in Type-I, we can also omit the case where $\gamma_{js} = 1$ for some vertex $j$. Contributions from all the $m_2$ Type-II components to $\psi_G$ are then bounded by

$$\left| \mathbb{E} \prod_{C_s \in \mathcal{F}_2} \psi_{C_s} \right| \leq K(\varepsilon_n n^{1/4})^{\sum_{s=1}^{m_2} \sum_{j=1}^{u_s} (\gamma_{js} - 2)}. \tag{7.14}$$

Combining (7.13) and (7.14) by noticing the relation $\sum_i \mu_i + \sum_{j,s} \gamma_{js} = 2k$, the overall contribution from random variables has a bound

$$|\mathbb{E} \psi_G| \leq K(\varepsilon_n n^{1/4})^{\sum_{i=1}^{m_1} (\mu_i - 4) + \sum_{s=1}^{m_2} \sum_{j=1}^{u_s} (\gamma_{js} - 2)} = K(\varepsilon_n n^{1/4})^{2k - 4m_1 - 2\sum_{s=1}^{m_2} u_s}. \tag{7.15}$$

Next the estimation of the weight part $a_G$ will be established. Since $\sum_{j=1}^{n} |a_{jj}|^w = \mathrm{O}(n)$ holds for any positive integer $w$, thus

$$\sum_{s_1, \ldots, s_{m_1}} \prod_{C \in \mathcal{F}_1} |a_C| = \sum_{s_1, \ldots, s_{m_1}} \prod_{i=1}^{m_1} |a_{s_i s_i}|^{\mu_i/2} \leq \prod_{i=1}^{m_1} \left( \sum_{k=1}^{n} |a_{kk}|^{\mu_i/2} \right) \leq K n^{m_1}. \tag{7.16}$$

For a given Type-II component $C_s$ with $t_s$ edges, $e_1, \ldots, e_{t_s}$ and $u_s$ vertices, $v_1, \ldots, v_{u_s}$, we extract a spanning tree from $C_s$ and assume its edges are $e_1, \ldots, e_{u_s - 1}$, without loss of generality. However, we need to distinguish two situations for the remaining sub-graph after extraction of the spanning tree. Let $\rho(A_n)$ be the spectral norm of the matrix $A_n$.

*Case* 1. The remaining sub-graph has at least one edge $u \to v$ with $u \neq v$. Note that

$$\sum_{v_1} |a_{v_1 v_2}|^2 \leq \rho(A_n)^2. \tag{7.17}$$

This, via induction, implies that we have for the tree part

$$\sum_{v_1,\ldots,v_{u_s}} \prod_{j=1}^{u_s-1} |a_{e_j}|^2 \leq \rho(A_n)^{2u_s-2} n,$$

and for the remaining sub-graph

$$\sum_{v_1,\ldots,v_{u_s}} \prod_{j=u_s}^{t_s} |a_{e_j}|^2 \leq \rho(A_n)^{2t_s-2u_s+2} n^{u_s-1}.$$

In the second inequality above, we use the fact that $t_s > u_s$ as all degrees of vertex of Type-II are no less than 2. It follows that

$$\sum_{v_1,\ldots,v_{u_s}} \prod_{j=1}^{t_s} |a_{e_j}| \leq \left( \sum_{v_1,\ldots,v_{u_s}} \prod_{j=1}^{u_s-1} |a_{e_j}|^2 \sum_{v_1,\ldots,v_{u_s}} \prod_{j=u_s}^{t_s} |a_{e_j}|^2 \right)^{1/2} \leq \rho(A_n)^{t_s} n^{u_s/2}, \tag{7.18}$$

which gives, together with (7.16), that

$$\sum_G |a_G| = \sum_G \prod_{C \in \mathcal{F}_1} |a_C| \prod_{C_s \in \mathcal{F}_2} |a_{C_s}| \leq K n^{u_1/2 + \cdots + u_{m_2}/2 + m_1}. \tag{7.19}$$

Combining (7.15) and (7.19), we obtain

$$n^{-k/2} \left| \mathbb{E} \sum_G a_G \psi_G \right| \leq n^{-k/2} \sum_G |a_G| |\mathbb{E}\psi_G| \leq K n^{-k/2} (\varepsilon_n n^{1/4})^{2k-4m_1-2\sum_{s=1}^{m_2} u_s} \sum_G |a_G|$$

$$\leq K \varepsilon_n^{2k-4m_1-2\sum_{s=1}^{m_2} u_s}. \tag{7.20}$$

*Case* 2. The remaining sub-graph does not contain any edge $u \to v$ with $u \neq v$, i.e., all remaining edges are loops. Since the degree of each vertex of a Type-II component is no less than two, there must exist at least two vertices whose degrees are more than two. Thus (7.14) turns into

$$\left| \mathbb{E} \prod_{C_s \in \mathcal{F}_2} \psi_{C_s} \right| \leq K (\varepsilon_n n^{1/4})^{\sum_{s=1}^{m_2} \sum_{j=1}^{u_s} (\gamma_{js}-2)-2m_2}. \tag{7.21}$$

We now need to consider two possibilities.

(a) If all vertices of a connected sub-graph have a loop, then similar to (7.18), we have

$$\sum_{v_1,\ldots,v_{u_s}} \prod_{j=1}^{t_s} |a_{e_j}| \leq \left( \sum_{v_1,\ldots,v_{u_s}} \prod_{j=1}^{u_s-1} |a_{e_j}|^2 \sum_{v_1,\ldots,v_{u_s}} \prod_{j=u_s}^{t_s} |a_{e_j}|^2 \right)^{1/2} \leq \rho(A_n)^{t_s} n^{(u_s+1)/2}$$

and then (7.19) becomes

$$\sum_G |a_G| \leq K n^{(u_1+1)/2 + \cdots + (u_{m_2}+1)/2 + m_1}. \tag{7.22}$$

But, at this point, there must exist a vertex such that its degree exceeds three and so (7.21), correspondingly, changes into

$$\left| \mathbb{E} \prod_{C \in \mathcal{F}_2} \psi_C \right| \leq K(\varepsilon_n n^{1/4})^{\sum_{s=1}^{m_2} \sum_{j=1}^{u_s} (\gamma_{j_s}-2)-4m_2}. \tag{7.23}$$

By (7.22) and (7.23), similar to (7.20), we get

$$n^{-k/2} \left| \mathbb{E} \sum_G a_G \psi_G \right| \leq Kn^{-m_2/2} \varepsilon_n^{2k-4m_1-2\sum_{s=1}^{m_2} u_s - 4m_2} \leq Kn^{-m_2/2}. \tag{7.24}$$

The last inequality results from the fact that by construction, the exponent of $\varepsilon_n$ is nonnegative. Consequently, the contributions from such graphs can be neglected.

(b) If not all vertices of a connected sub-graph have a loop, then

$$\sum_{v_1,\dots,v_{u_s}} \prod_{j=1}^{t_s} |a_{e_j}| \leq \rho(A_n)^{t_s} n^{u_s/2},$$

and, correspondingly, (7.19) becomes

$$\sum_G |a_G| \leq Kn^{u_1/2 + \cdots + u_{m_2}/2 + m_1}. \tag{7.25}$$

To see it, as an example, we consider the following

$$\left| \sum_{v_1,v_2,v_3} a_{v_1 v_2} a_{v_1 v_3} a_{v_2 v_2} a_{v_3 v_3} \right| = \left| \sum_{v_2,v_3} b_{v_2 v_3} a_{v_2 v_2} a_{v_3 v_3} \right| \leq \left( \sum_{v_2,v_3} |b_{v_2 v_3}|^2 \right)^{1/2} \left( \sum_{v_3} a_{v_3 v_3}^2 \right)^{1/2} \left( \sum_{v_2} a_{v_2 v_2}^2 \right)^{1/2}$$
$$= O(n^{3/2}),$$

where

$$b_{v_2 v_3} = \sum_{v_1} a_{v_1 v_2} a_{v_1 v_3} = \sum_{v_1} \bar{a}_{v_2 v_1} a_{v_1 v_3} = (\overline{A_n} A_n)_{v_2 v_3}.$$

For general cases, we can verify the order by induction. Using (7.13), (7.21) and (7.25), similar to (7.20), we get

$$n^{-k/2} \left| \mathbb{E} \sum_G a_G \psi_G \right| \leq Kn^{-m_2/2} \varepsilon_n^{2k-4m_1-2\sum_{s=1}^{m_2} u_s - 2m_2} \leq Kn^{-m_2/2}. \tag{7.26}$$

So the contribution from this kind of graph can also be neglected. Here we remind the reader that (7.20) is obtained by assuming all $m_2$ Type-II components belonging to case 1 and that (7.24) or (7.26) holds if all $m_2$ Type-II components belong to case 2. If some Type-II components of the graph $G$ belong to case 1 and the others pertain to case 2, by a similar argument to the above, one can show that

$$n^{-k/2} \left| \mathbb{E} \sum_G a_G \psi_G \right| = o(1). \tag{7.27}$$

Therefore, if some item involves the connected sub-graph of case 2, the contribution from this item can then be omitted.

In summary, in conjunction with (7.20) and the meanings of $2k, 4m_1, 2\sum_{s=1}^{m_2} u_s$, we know that the graphs leading to a non negligible term are as follows: then degrees of vertices of all its Type-I components must be four; its Type-II components all fall into case 1 such that all its vertices are of degree two. More precisely, we know that only the following isomorphic classes give a dominating term:

- $k_1$ double loops $u \to u$ with terms $a_{uu}^2[\sum_{\ell=1}^K c_\ell(\bar{x}_{\ell u}y_{\ell u} - \rho(\ell))]^2$;
- $k_2$ simple cycles $u \to v, v \to u$ with terms $|a_{uv}|^2[\sum_{\ell=1}^K \bar{x}_{\ell u}y_{\ell v}][\sum_{\ell=1}^K \bar{x}_{\ell v}y_{\ell u}]$;
- $k_3$ double arrows $u \to v, u \to v$ with terms $a_{uv}^2[\sum_{\ell=1}^K c_\ell \bar{x}_{\ell u}y_{\ell v}]^2$.

In addition, the degrees of vertices satisfy

$$4(k_1 + k_2 + k_3) = 2k,$$

which implies that $k$ must be even. Therefore, let $k = 2p$ be an even integer. We notice that here, the relations on the edges, namely $2(k_1 + k_2 + k_3) = k$, hold automatically. Thus, we can claim that

$$\mathbb{E}\xi_n^{2p} = \frac{1}{n^p} \sum_{k_1+k_2+k_3=k} \frac{(2p)!}{k_1!k_2!k_3!} \times C_1 \times C_2 \times C_3 + \mathrm{o}(1), \tag{7.28}$$

where

$$C_1 = \prod_{j=1}^{k_1} \mathbb{E}\left[a_{u_j u_j}^2 \left\{\sum_{\ell=1}^K c_\ell(\bar{x}_{\ell u_j}y_{\ell u_j} - \rho(\ell))\right\}^2\right], \qquad \{u_j\} \subset \{1,2,\ldots,n\},$$

$$C_2 = \prod_{j=1}^{k_2} \mathbb{E}\left[|a_{u_j v_j}|^2 \left\{\sum_{\ell=1}^K c_\ell(\bar{x}_{\ell u_j}y_{\ell v_j})\right\}\left\{\sum_{\ell=1}^K c_\ell(\bar{x}_{\ell v_j}y_{\ell u_j})\right\}\right], \quad \{u_j,v_j\} \subset \{1,2,\ldots,n\},$$

$$C_3 = \prod_{j=1}^{k_3} \mathbb{E}\left[a_{u_j v_j}^2 \left\{\sum_{\ell=1}^K c_\ell \bar{x}_{\ell u_j}y_{\ell v_j}\right\}^2\right], \qquad \{u_j,v_j\} \subset \{1,2,\ldots,n\}.$$

Let

$$\alpha_1 = \mathbb{E}\left[\sum_{\ell=1}^K c_\ell(\bar{x}_{\ell 1}y_{\ell 1} - \rho(\ell))\right]^2, \tag{7.29}$$

$$\alpha_2 = \mathbb{E}\left[\left\{\sum_{\ell=1}^K c_\ell(\bar{x}_{\ell 1}y_{\ell 2})\right\}\left\{\sum_{\ell=1}^K c_\ell(\bar{x}_{\ell 2}y_{\ell 1})\right\}\right], \tag{7.30}$$

$$\alpha_3 = \mathbb{E}\left[\sum_{\ell=1}^K c_\ell \bar{x}_{\ell 1}y_{\ell 2}\right]^2. \tag{7.31}$$

By (7.20) or (7.27) again, along with inclusion-exclusion principle, (7.28) turns into

$$\mathbb{E}\xi_n^{2p} = \frac{1}{n^p}(2p-1)!!$$

$$\times \sum_{k_1+k_2+k_3=k} \frac{(p)!}{k_1!k_2!k_3!} \prod_{(j_1,j_2,j_3)=(1,1,1)}^{(k_1,k_2,k_3)} \alpha_1^{k_1}\alpha_2^{k_2}\alpha_3^{k_3} a_{u_{j_1}u_{j_1}}^2 |a_{u_{j_2}v_{j_2}}|^2 a_{u_{j_3}v_{j_3}}^2 + \mathrm{o}(1)$$

$$= \frac{1}{n^p}(2p-1)!!\left(\alpha_1 \sum_{u=1}^n a_{uu}^2 + \alpha_2 \sum_{u \neq v} |a_{uv}|^2 + \alpha_3 \sum_{u \neq v} a_{uv}^2\right)^p + \mathrm{o}(1).$$

### 7.2. Truncation

The truncation and renormalization under the fourth-moment condition is by now standard, see e.g. [3]. For our purpose and for ease of presentation, we give full details in the case of $K = 1$. The general case goes through in a same manner.

We aim to the replacement of the entries of $x$ and $y$ with truncated, centralized, normalized variables. Let $\hat{x} = (\hat{x}_1, \ldots, \hat{x}_n)^{\mathrm{T}}$ and $\tilde{x} = (\tilde{x}_1, \ldots, \tilde{x}_n)^{\mathrm{T}}$, where

$$\hat{x}_j = x_j I(|x_i| \leq \varepsilon_n n^{1/4}), \qquad \tilde{x}_j = \hat{x}_j - \mathbb{E}\hat{x}_j, \quad j = 1, \ldots, n.$$

Since $\mathbb{E}|x_1|^4 < \infty$, for any $\varepsilon > 0$

$$n\mathbb{P}(|x_1| \geq \varepsilon n^{1/4}) \to 0,$$

and then, because of the arbitrariness of $\varepsilon$, there exists a positive sequence $\varepsilon_n$ such that

$$n\mathbb{P}(|x_1| \geq \varepsilon_n n^{1/4}) \to 0 \quad \text{and} \quad \varepsilon_n \to 0.$$

It follows that

$$\mathbb{P}(x^* A_n y \neq \hat{x}^* A_n y) \leq n\mathbb{P}(|x_1| \geq \varepsilon_n n^{1/4}) = \mathrm{o}(1). \tag{7.32}$$

For $h = 1, 2, \ldots$, find $n_h (n_h > n_{h-1})$, for all $n > n_h$ with

$$h^{12} \int_{|x_1| \geq \sqrt[4]{n}/h} |x_1|^4 < 2^{-h}.$$

Let $\rho_n = \frac{1}{h}$ for all $n \in [n_h, n_{h+1}]$, thus, as $n \to \infty$, $\rho_n \to 0$ and

$$\rho_n^{-12} \int_{|x_1| \geq \sqrt[4]{n}\rho_n} |x_1|^4 \to 0. \tag{7.33}$$

Now, for each $n$, let $\gamma_n$ be the larger of $\rho_n$ and $\varepsilon_n$. However, in the following we still use the notation $\varepsilon_n$ instead of $\gamma_n$. By Markov inequality and Burkholder inequality

$$\begin{aligned}
\mathbb{P}(|\hat{x}^* A_n y - \tilde{x}^* A_n y| \geq \delta) &\leq |\mathbb{E}[x_1 I(|x_1| \leq \varepsilon_n n^{1/4})]|^4 \frac{\mathbb{E}|\sum_{i=1}^n y_i(\sum_{j=1}^n a_{ji})|^4}{\delta^4} \\
&\leq K\mathbb{E}[|x_1|^4 I(|x_1| \geq \varepsilon_n n^{1/4})]n^{-2}\varepsilon_n^{-12} \frac{[\sum_{i=1}^n |\sum_{j=1}^n a_{ji}|^2]^2}{\delta^4} \\
&\leq K\mathbb{E}|x_1|^4 I(|x_1| \geq \varepsilon_n n^{1/4})\varepsilon_n^{-12},
\end{aligned} \tag{7.34}$$

where we have used the inequality

$$\sum_{i=1}^n \left|\sum_{j=1}^n a_{ji}\right|^4 \leq \left[\sum_{i=1}^n \left|\sum_{j=1}^n a_{ji}\right|^2\right]^2$$

and the fact

$$\sum_{i=1}^n \left|\sum_{j=1}^n a_{ji}\right|^2 = \mathrm{O}(n).$$

From (7.34) and (7.33) we have

$$\hat{x}^* A_n y - \tilde{x}^* A_n y \xrightarrow{\text{i.p.}} 0. \tag{7.35}$$

Next, we need to normalize the truncated variables $\tilde{x}_i$'s. It is evident that $\lim_{n\to\infty} \mathbb{E}|\tilde{x}_1|^2 = 1$ and that

$$\begin{aligned}
|1 - \sqrt{\mathbb{E}|\tilde{x}_1|^2}| &\leq |1 - \mathbb{E}|\tilde{x}_1|^2| \leq 2\mathbb{E}[|x_1|^2 I(|x_1| \geq \varepsilon_n \sqrt[4]{n})] \\
&\leq 2\varepsilon_n^{-2} n^{-1/2}\mathbb{E}[|x_1|^4 I(|x_1| \geq \varepsilon_n \sqrt[4]{n})],
\end{aligned} \tag{7.36}$$

which, together with (7.33), gives that

$$\left| (1 - \sqrt{\mathbb{E}|\tilde{x}_1|^2}) \frac{\operatorname{tr} A_n}{\sqrt{n}} \right| \le 2\varepsilon_n^{-2} \mathbb{E}|x_1|^4 I(|x_1| \ge \varepsilon_n \sqrt[4]{n}) \to 0. \tag{7.37}$$

Combining (7.32), (7.35) and (7.37), it is now sufficient to consider

$$\frac{1}{\sqrt{n}} \left( \frac{\tilde{x}^* A_n y}{\sqrt{\mathbb{E}|\tilde{x}_1|^2}} - \rho \operatorname{tr} A_n \right)$$

instead of $\frac{1}{\sqrt{n}}(x^* A_n y - \rho \operatorname{tr} A_n)$. Moreover, it is not difficult to see that we can substitute $\rho' = \operatorname{cov}(\tilde{x}_1 / \sqrt{\mathbb{E}|\tilde{x}_1|^2}, y_1)$ for $\rho$ without alternating the weak limit.

The truncation, centralization and normalization of $y$ can be completed with a similar argument as above. In the sequel, for simplicity, we shall suppress all superscripts on the variables and suppose that $|x_i| \le \varepsilon_n n^{1/4}$, $|y_i| \le \varepsilon_n n^{1/4}$, $\mathbb{E}x_i = \mathbb{E}y_i = 0$, $\mathbb{E}|x_i|^2 = \mathbb{E}|y_i|^2 = 1$, $\mathbb{E}|x_1|^4 < \infty$, $\mathbb{E}|y_1|^4 < \infty$ and we still denote by $\rho$ the covariance between the transformed variables.